\documentclass[11pt,a4paper]{amsart}
\usepackage{setspace}
\usepackage{amsfonts}
\usepackage{amssymb}
\usepackage{amsmath}
\usepackage{yhmath}
\usepackage{xy}
\usepackage{enumerate}
\usepackage{amsthm}
\usepackage[mathscr]{eucal}

\numberwithin{equation}{subsection}

\begin{document}
\newtheorem{thm}{Theorem}[subsection]
\newtheorem{tprop}[thm]{Proposition}
\newtheorem{tlemma}[thm]{Lemma}
\newtheorem{tcor}[thm]{Corollary}

\newtheorem{prop}{Proposition}[subsection]
\newtheorem{pthm}[prop]{Theorem}
\newtheorem{plemma}[prop]{Lemma}
\newtheorem{pcor}[prop]{Corollary}

\newtheorem{lemma}{Lemma}[subsection]
\newtheorem{lthm}[lemma]{Theorem}
\newtheorem{lprop}[lemma]{Proposition}
\newtheorem{lcor}[lemma]{Corollary}

\input xy
\xyoption{all}
\title{The Plancherel formula of $L^2(N_0 \setminus G;\psi)$}
\author{Tang U-Liang}
\address{Department of Mathematics,
National University of Singapore,
Block S17,
10, Lower Kent Ridge Road,
Singapore 119076 }
\email{mattul@nus.edu.sg}
\begin{abstract}
We study the right regular representation on the space $L^2(N_0\setminus G;\psi)$ where $G$ is a quasi-split $p$-adic group and $\psi$ a non-degenerate unitary character of the unipotent subgroup $N_0$ of a minimal parabolic subgroup of $G$. We obtain the direct integral decomposition of this space into its constituent representations. In particular, we deduce that the discrete spectrum of $L^2(N_0\setminus G;\psi)$ consists precisely of $\psi$ generic discrete series representations and derive the Plancherel formula for $L^2(N_0\setminus G;\psi)$. 
\end{abstract}
\keywords{Whittaker model, generic representations, square integrable representations}
\thanks{The author wishes to thank Professor Gordon Savin, Professor Gan Wee Teck  and H.Y. Loke for initiating this project, suggestions, comments and guidance. }
\maketitle


\begin{raggedright}
\parindent=0.5in

\setcounter{tocdepth}{1}
\tableofcontents

\section{Introduction and statement of main results}
\renewcommand{\hom}{\mathrm{Hom}}
\newcommand{\im}{\mathrm{Im}}
\newcommand{\ur}{\mathrm{ur}}
\newcommand{\ra}{\rangle}
\newcommand{\la}{\langle}
\renewcommand{\v}[1]{v_{\nu}(#1)}
\newcommand{\w}[1]{w_{\nu}(#1)}
\newcommand{\intorbit}{\int_{\mathcal{O}}}
\newcommand{\deltap}{\delta_{\bar{P}_{\theta}}}
\newcommand{\unnind}{\mathrm{Ind}}

Let $\psi$ be a nondegenerate  unitary character of the unipotent radical $N_0$ of a minimal standard parabolic subgroup of a connected quasi-split $p$-adic group, $G$. Define $L^2(N_0\setminus G;\psi)$ as the space of functions on $G$ which  transform according to $\psi$, i.e. $f(ng)=\psi(n)f(g)$ and are square integrable modulo $N_0$. This space becomes a unitary representation of $G$ via right translation.   

The purpose of this work is to obtain the Plancherel formula for this unitary representation. Dinakar Ramakrishnan has obtained the spectral decomposition for this space when $G=GL(2)$ in \cite{ram2}. He proves this for both the archimedean and non-archimedean cases. Nolan Wallach then extended his result to arbitrary real reductive groups. Indeed, much of the arguments in the latter part of this paper imitate Wallach's original proofs in \cite{wallach}. 

There are two crucial components into developing a spectral theory of $L^2(N_0 \setminus G; \psi)$. One of them is to prove a certain surjectivity result from the space of Schwartz functions $G$ to the space of Schwartz functions on $N_0 \setminus G$. This is the author's original contribution. 

A second component needed is the knowledge of the analyticity of a certain Jacquet integral (this integral is defined later in the paper.) In the non-archimedean case, it follows from an argument of Casselman in \cite{casselman} (see also \cite{jacq} and \cite{shah}) that the Jacquet integral extends to a holomorphic function. One requires such an integral to define a \emph{Whittaker transform} which transforms certain smooth functions on orbits of discrete series representations to Schwartz functions on $N_0 \setminus G$  analogous to the Harish-Chandra wave packet map. We refer the reader to  Section \ref{sec62} for the precise definition of this transform. 

It is a pleasant consequence of our investigations that we are able to obtain an extension of a result in \cite{sav}. We refer the reader to Theorem \ref{main2}. This generalization is a consequence of the fact that matrix coefficients of generic discrete series representations decay rapidly (in manner made precise in this paper). 

To state our main result, we must describe briefly the notation involved.  Let $P=MN$ be a standard parabolic subgroup of $G$ with $M$ and $N$ its Levi and unipotent subgroup respectively Let $\psi^{M}$ denote the restriction of $\psi$ to $M \cap N_0$. We take a $\psi^M$ generic discrete series representation $(\sigma, H_{\sigma})$ and consider the unitarily induced representation $I_P^G(\sigma\otimes \nu)$ where $\nu \in \mathrm{Im}(X_{ur}(M))=i\mathfrak{a}_M/L$ ($L$ is some lattice in the real vector space $\mathfrak{a}_M$) runs over all unramified unitary  characters of $M$. Let $Wh_{\psi^M}(H_{\sigma})$ denote the (one dimensional) space of Whittaker functionals on $\sigma$. Let $$\mathscr{H}_{\sigma, \nu}=I_P^G(\sigma, \nu) \otimes Wh_{\psi^M}(H_{\sigma})$$ and consider the direct integral $$\mathcal{I}_{\sigma, M}=\int^{\oplus}_{i\mathfrak{a}_M/L} \mathscr{H}_{\sigma, \nu}\, \tilde{\mu}(\sigma, \nu)\,d\nu$$ where $\tilde{\mu}(\sigma, \nu)$ is a certain normalization of the Plancherel measure on $i\mathfrak{a}/L$. 

Let $$W(G|M):=\{w \in W^G\mid w.M=M\}/W^M$$ where $W^G$ and $W^M$ denote the Weyl group of $G$ and its Levi $M$ respectively. If $\sigma \in \mathscr{E}^2_{\psi^M}$, then $w.\sigma$ is defined and denote $\mathscr{E}^2_{\psi^M}(M)/W(G|M)$ to be the set of isomorphism classes of square integrable representations of $M$ which are $\psi^M$ generic modulo the action of $W(G|M)$. 

We prove that 
\begin{thm} \label{main}
There exists a unitary linear surjection from $$\sum_{M \subset G}\sum_{\sigma \in\mathscr{E}^2_{\psi^M}(M)/W(G|M)} \mathcal{I}_{\sigma,M} $$ onto $L^2(N_0 \setminus G; \psi).$ 
\end{thm}

We also compute the explicit normalization of $\tilde{\mu}(\sigma, \nu)$. If $\mu(\sigma,\nu)\,d\nu$ denotes the Plancherel measure on $i\mathfrak{a}_M/L$, then $$\tilde{\mu}(\sigma,\nu)\,d\nu = \frac{1}{|W(G|M)| \gamma(G|M)c^2(G|M)}\mu(\sigma,\nu)\,d\nu.$$

This paper is organized as follows. In Section 2 we give a simplified exposition of Bruhat-Tits theory adequate for our purposes and prove  Lemma \ref{nomatchingcoset}. This is a key lemma required to prove the crucial Lemma \ref{pairing}. 

Section 3 describes the notations and conventions for parabolic subgroups and tori needed to describe unitary parabolic induction. We also give a description of both the Schwartz spaces on $G$ and $N_0 \setminus G$ needed later. 

In Section 4 we discuss the asymptotics of Whittaker functions for discrete series representations and more generally for tempered representations. 
 
Finally in Section 5 we prove Lemma \ref{pairing} which tells us (among other things) that the map from $\mathcal{C}^*(G)$ to $\mathcal{C}^*(N_0 \setminus G; \psi)$ is surjective.  

In Section 6, we set the stage for deriving the Plancherel formula. This is essentially Theorem \ref{main4}. Our main result, Theorem \ref{main}, is precisely Corollary \ref{lastcor}. 

\textbf{Remark}. While this paper was being written, it was brought to the author's attention that Erez Lapid and Mao Zhengyu  had obtained an explicit form of the Whittaker function and its asymptotics on a split group $G$ in \cite{lapid}. Theorem \ref{suff} is a direct corollary of their results.

They conjectured the following: Let $W(\pi)$ denote the Whittaker model of a generic representation and suppose that $\int_{Z_G N_0\setminus G} |W(g)|^2 \;dg$ is finite for all $W(g)\in W(\pi)$,  then $\pi$ is , square integrable.  By Theorem \ref{location of image} and Theorem \ref{main3} one concludes that this conjecture is true. 

Patrick Delorme has obtained the results of this work independantly in \cite{delorme} and \cite{delormepaleyweiner}. However, our approach differs slightly from his treatment. We also thank Professor Delorme for pointing out a gap in the previous version of Proposition \ref{corrected}.

\section{The Cartan and Iwasawa decompositions of $G$}

\newcommand{\val}{\mathrm{val}}
\newcommand{\ncoset}{[N_0/ N_0 \cap g_i Hg_i^{-1}]}
\newcommand{\schwartzg}{\mathcal{C}^*(G)^H}
\subsection{}\label{intro51}
We begin by fixing a $p$-adic field $k$ with ring of integers $\mathfrak{o}$ and normalized absolute value $| .| =q^{-\val(.)}$.
Let $\mathbf{G}$ be a connected quasi split reductive group  defined over $k$. Let $G=G(k)$ be its $k$-rational points.
Let $A_0$ the maximal $k$ split torus and $M_0$ its (abelian) centralizer with $M_{0,0}$ the maximal compact subgroup of $M_0$. 

With $W$ denoting the Weyl group of $G$ with respect to $A_0$, let  $\widetilde{W}$ be the affine Weyl group extending the Weyl group $W$. We may identify $\widetilde{W}$ as the semidirect product of $W$ and $D = A_0/A_0^1$ where $A_0^1$ is the maximal open compact subgroup of $A_0$. Let $\Sigma_{\mathrm{aff}}$ and $\,^{nd}\Sigma$ denote the affine root system and set of nondivisible roots of $G$ respectively. Let $K$ be the special maximal open compact subgroup fixing a special point $x_0$ in the apartment stabilized by $A_0$. Denote by $\Sigma_0$ the roots of $\Sigma_{\mathrm{aff}}$ vanishing on $x_0$. If $B$ is the subgroup of $G$ fixing a chamber pointwise in the apartment with vertex $x_0$, then $B$ is an Iwahori subgroup. One knows that there is a bijection $\lambda:\,^{nd}\Sigma\rightarrow \Sigma_0$ since every root $\alpha \in \Sigma_0$ is a \emph{positive} multiple of a unique root in $\,^{nd}\Sigma$. We write $\lambda(\alpha)$ as $\lambda_{\alpha}\alpha$. 

The Iwahori decomposition for $G$ and $K$ is $$G=\coprod_{w \in \widetilde{W}} B w B$$ and $$K=\coprod_{w\in W} BwB$$ respectively with both these unions disjoint.

We will write the Iwasawa decomposition of $G$ as $N_0A_0 K$ where $N_0$ is the unipotent subgroup of a choice of minimal parabolic $P_0$ of $G$ containing $M_0$. 

Define $\,^{nd}\Sigma^+$ to be the system of positive roots determined by $P_0$ and let $\Delta$ denote the set of simple roots.  For $\epsilon > 0$, define  $$A^+(\epsilon^{-1}):=\{a \in A_0 \mid |\alpha(a)| \leq \epsilon^{-1} \, \forall \alpha  \in \Delta\}.$$ Writing $A^{+}$ as $A^+(1)$, we have the Cartan decomposition $G=KA^+K$.

Write $\{H_m\}_{m\geq1}$ for the system of \lq good\rq\  open compact subgroups for the topology on $G$ (see \cite{sil}). Without loss of generality we assume that if $w\in W$, then $w \notin H_i$ for any $i \geq1$. 

For each $\alpha \in \Sigma_{\mathrm{aff}}$ define $N_0(\alpha)$ as in \cite{casselunram}. Then for $\alpha \in \,^{nd}\Sigma^+$ and $m \in \mathbb{Z}$ define $$N_{\alpha, m} := N_0(\lambda(\alpha)+m).$$ We have that $N_{\alpha,m+1}\subsetneq N_{\alpha, m}$ and if $N_{\alpha}:= \bigcup_{m \in \mathbb{Z}} N_{\alpha, m}$, then $$N_0= \prod_{\alpha \in \,^{nd} \Sigma^+}\,N_{\alpha}$$ in any order.

 Write $N_{0,m}^-=\prod_{\alpha \in \,^{nd}\Sigma^+}N_{-\alpha,m}$, $N_{0,m}=\prod_{\alpha \in \,^{nd}\Sigma^+}N_{\alpha,m}$. Then $B=N^-_{0,1}M_{0,0} N_{0,0}.$ If $m\geq 0$, $H_m = N^-_{0,m} (M_0 \cap H_m) N_{0,m}$. 
 
Set $q_{\alpha}=[N_0(\alpha-1) : N_0(\alpha)]$ and $q_{\alpha/2}=q_{\alpha+1}/q_{\alpha}$. We remark that it is possible that $q_{\alpha/2} \neq 1$. For convenience,  $q_{\alpha/2}^{\lfloor m/2 \rfloor}q_{\alpha}^m=[N(\alpha +1): N(\alpha+m+1)]$ is shortened to $\,_{\alpha}q^m.$

Fix an open compact subgroup $H=H_m$ as defined above of $K$. We determine a parametrization of $G/H$ cosets in a single $N_0\setminus G/H$ coset. Let $\{g_i\}$ be a set of coset representatives of $N_0 \setminus G/H$. By the Iwasawa decomposition we may assume $g_i=ak$ where $a \in A$ and $k \in K$ where $k$ comes from a set of coset representatives of $K/H$. 

Under the left action of $N_0$, the stabilizer of $g_iH$ is $N_0 \cap g_iHg_i^{-1}$. Let $[N_0/ N_0 \cap g_i Hg_i^{-1}]$ be a set of $N_0 \cap g_i Hg_i^{-1}$-coset representatives in $N_0$. Then each $G/H$ coset in $N_0g_iH$ is parametrized the set $[N_0/ N_0 \cap g_i Hg_i^{-1}]$.

Since $H$ is normal in $K$, $N_0 \cap g_iHg_i^{-1}= N_0 \cap aHa^{-1}=aN_{0,m}a^{-1}.$ Written component wise $$aN_{\alpha,m}a^{-1}=N_0(\lambda(\alpha)+m+\lambda_{\alpha}\mathrm{val}(\alpha(a)))$$ where we recall that $\lambda_{\alpha}$ is the positive multiple arising from the bijection between $\,^{nd}\Sigma$ and $\Sigma_0$. 

Now chose an integer $n$ so that $N_{0,n} \supsetneq N_{0,m}$. Also, choose $\epsilon>0$ small enough so that if $a$ satisfies $$\mathrm{val}(\alpha(a)) \geq \frac{1}{\lambda_{\alpha}}(-m+n)$$  for all simple roots $\alpha \in \Delta$ then $a \in A^+(\epsilon^{-1})$.  Then $a \in A^+(\epsilon^{-1})$ contains the set of all $a \in A$ such that $N_0\cap aHa^{-1} \subset N_{0,n}$. 

\begin{lemma}\label{nomatchingcoset}
Assume we are given a fixed $a \in A^+(\epsilon^{-1})$ and fixed $k \in K.$  Then  for any $n \in N_0- ((N_0\cap aH_ma^{-1})N_{0,m})$, $$H_mnakH_m \neq H_mwa'w^{-1}kH_m$$ for any $a'\in A^+$ and any $w\in W$. 
\end{lemma}
\begin{proof}
Firstly we have to compute the $KA^+K$ decomposition for the element $g=nak$ as we run over $n \in N_0$ for a fixed $ak$.

Fix an $a \in A^+(\epsilon^{-1})$ such that $N_0\cap aH_ma^{-1} \nsubseteq N_{0,m}$. There are three cases. 

\begin{enumerate} 
\item If $n \in N_{0,0}-((N_0\cap aH_ma^{-1})N_{0,m})$, then $$g=nak=(nw)a^+(w^{-1}k)$$ where $nw$ and $w^{-1}k$ are in $K$. 
\item  If $n \notin N_{0,0}$ but $a^{-1}na^{} \in K$, then $nak= a (a^{-1} n a) k = wa^+w^{-1}n_0 k$ where \newline $n_0 = a^{-1} na \in K$. As $w^{-1} n_0k \in K$, the $KA^+K$ decomposition is $$nak= (w) a^{+} (w^{-1} n_0k).$$ 
\item If $n \in N_0$ is not of the either two cases then we consider the Iwahori decomposition of $(na)$ so that $(na)k=(b_{na}\tilde{w}b'_{na})k$  where $\tilde{w} \in \widetilde{W}$ and $b_{na}, b'_{na} \in B$. Since $\widetilde{W}=W.D$, $(na)k=(b_{na} w' a_1 b'_{na})k$ where  we write $\tilde{w}=w'a_1$ with $w'\in W$ and $a_1 \in D$. Since $na$ and $a^{-1}na \notin K$, $w'\neq 1$ and $a_1\neq 1$. 
Then $$nak=(b_{na}w'w_1)a_1^{+}(w_1^{-1}b'_{na}k)$$ where $a_1^{+}=w_1^{-1}a_1w_1$ with $w_1$ chosen so that $a_1^{+} \in A^+$.  The bracketed elements are elements in $K$. 
\end{enumerate}

For case (1), let us assume that there exists $w_0 \in W$ and $a_0 \in A^+$ such that $$(nw)a^+(w^{-1}k) \in H_m w_0 a_0 w_0^{-1}kH_m.$$  
Since the Cartan decomposition gives a disjoint union over $K\check{a}K$ where $\check{a}$ runs over $[A^+/A_0^1]$ $a^+=a_0$ and thus $nw \in H_mw_0$. But this means that $nww_0^{-1} \in H_m$ which is impossible since $n \notin H_m$. This proves the lemma in this case. 

For case (3) we argue in a similar fashion. So suppose $$(b_{na}w'w_1) a_1^+ (w_1^{-1} b'_{na} k) \in H_m w_0a_0 w_0^{-1} kH_m.$$ As before, $a_1^{+}=a_0$ so that $b_1w'w_1 \in H_mw_0$. We may as well assume $b_1 \in H_m$ implying that $w'_1=w_0$. Now on the other hand $$w_1^{-1}b_{na}'k \in w_0^{-1} kH_m$$ implying that $$w_1^{-1} b_{na}' \in w_0^{-1}H_m.$$ Once again, we may assume  $b_{na}'\in H_m$ so that $w_0=w_1$. Thus $w'=1$ which is a contradiction.

Finally in order to apply the argument of the previous cases above to case (2), we must exclude the possibility that $n_0$ is in $H_m$. 

Indeed, if $n_0 \in H_m$, we will derive a contradiction. Write the component of an arbitrary $n\in N_0$ in $N_{\alpha}$ as $n_{\alpha}$. Since $a^{-1}na \in K$, $n_{\alpha} \in N_{\alpha}$ is contained in at most $N_{\alpha, \lambda_{\alpha} \mathrm{val}(\alpha(a))}$. However by our assumption on $n_0$, $n_{\alpha} \in N_{\alpha, \lambda_{\alpha}\mathrm{val}(\alpha(a))}$ satisfies $a^{-1}n_{\alpha}a \in N_{\alpha,m}$ for all roots $\alpha \in \,^{nd}\Sigma^+$. This implies that $n_{\alpha} \in N_{\alpha, \lambda_{\alpha} \mathrm{val}(\alpha(a))+m}$ i.e. $n \in N_0 \cap aHa^{-1}$. This contradicts the hypothesis of the lemma.

If $a \in A^+$ then $N_{0,m} \supset N_0\cap aHa^{-1}$  so that for any $n \notin N_{0,m}$, $nak$ is not in the same $H_m$-bicoset of as $wa'w^{-1}k$ for any $a'\in A^+$ and any $w\in W$. 
\end{proof}
\section{Parabolic subgroups, Tori and Schwartz spaces}
\subsection{}
We introduce the notion of standard parabolic subgroups and define a dense subspace of $L^2(N_0 \setminus G; \psi)$ which we will study in this paper. 

Fix a minimal parabolic subgroup $P_0$ as in the previous section. A \emph{standard parabolic pair} is a pair $(P, A)$ consisting of a parabolic subgroup $G \supset P \supset P_0$ and $A_0 \supset A \supset Z_G$ where $Z_G$ denotes the (split component of the) center of $G$. 

It is known that all such pairs are in one-to-one correspondence with subsets of $\Delta$. By abuse of notation we identify $(P,A)$ with $P$ and write this correspondence as $$\theta \mapsto P_{\theta}$$ where we agree that $P_{\emptyset}=P_0$ and $P_{\Delta}=G$. 
 It is well known that any standard parabolic corresponds to a subset $\theta\subset \Delta$. Conversely, to each subset $\theta$, one can associate a standard parabolic $P_{\theta}$, $M_{\theta}$, its Levi and
If $P_{\theta}$ is a standard parabolic, we write its Langlands decomposition as $P_{\theta}=M_{\theta}N_{\theta}$. Where there is no cause for confusion, we will drop the $\theta$ from notation. Moreover where the context is clear, $N_0$ always denotes $N_{\emptyset}$. 

A smooth unitary character $\psi$ of $N_0$ is said to be nondegenerate if and only if for any $\alpha \in \Delta$, $\psi$ restricted to $N_{\alpha}$ is non-trivial. 

\subsection{}\label{sec32}
Let $\delta_{P}$ denote the modular character of the parabolic subgroup $P$. By the Iwasawa decomposition $G = N_{\theta}M_{\theta}K$, we write $g \in G$ as $g = n_{P_{\theta}}(g)m_{P_{\theta}}(g)k(g)$. 

We define the constant $$\gamma(P_{\theta})=\int_{\bar{N_{\theta}}} \deltap(m_{P_{\theta}}(\bar{n})) d\bar{n}$$ 
where $\bar{N}_{\theta}$ denotes the unipotent subgroup of the opposite parabolic of $P_{\theta}.$ (c.f. \cite[pg. 240]{walsp}.) It is known that $\gamma(P_{\theta})$ does not depend on the choice of $N_{\theta}$ used to define it. Thus we may write $\gamma(P_{\theta})=\gamma(G|M_{\theta}).$ If $\alpha \in \,^{nd}\Sigma(P_{\theta}, A_{\theta})$, let $A_{\alpha}$ denote the (identity component) of the kernel of $\alpha$ and $M_{\alpha}$ the centralizer of $A_{\alpha}$ in $G$. Then $M_{\alpha} \supset M_{\theta}$ for every root $\alpha  \in \,^{nd}\Sigma(P_{\theta},A_{\theta})$.  Define
$$ c(G|M_{\theta})=\gamma(G|M_{\theta})^{-1} \prod_{\alpha \in \,^{nd}\Sigma(P_{\theta}, A_{\theta})} \gamma(M_{\alpha}|M_{\theta}).$$

If $(\pi,V)$ is a unitary representation, let $\la, \ra_{\pi}$ denote the Hermitian inner product on $V$. 

\begin{tlemma}\label{gammafact} If $(\sigma, H_{\sigma})$ is any unitary representation then $I_{\bar{P}_{\theta}}^G\sigma=I(\sigma)$ is a unitarily induced representation from $\bar{P}_{\theta}$ to $G$. Furthermore
for any $u,v \in I(\sigma)$ $$ \int_K \la u(k), v(k) \ra_{\sigma} \,dk= \gamma(G|M_{\theta})^{-1} \int_{N_{\theta}} \la u(n), v(n) \ra_{\sigma} \,dn. $$ 
\end{tlemma}
\begin{proof}
For the purposes of this proof,  $\unnind_P^G\; \sigma$ will denote unnormalized induction instead of the usual normalized induction. Recall the integral decomposition formula (c.f. \cite[pg. 240]{walsp}). 
$$\int_G f(g) \,dg = \gamma(G|M)^{-1}\int_{N \times M \times \bar{N}} f(nm\bar{n})\delta_P^{-1}(m)\,d\bar{n}\,dm\,dn$$ for any parabolic subgroup $P$ and $f \in C^{\infty}_c(G)$. 
 
For any $f \in C^{\infty}_{c}(G)$ let $P_{\delta_{\bar{P}_{\theta}}}(f)$ denote the projection of $f$ into $\unnind_{\bar{P}_{\theta}}^{G}\, \deltap$. This is given by $$\int_{\bar{P}_{\theta}} \deltap^{-1}(p) f(p)\, dp.$$ Let  $$I_{N_{\theta}}(w):= \int_{N_{\theta}} w(n)\,dn$$ for any $w\in \unnind_{\bar{P}_{\theta}}^{G} \deltap$ provided this integral converges.  

Consider \begin{align*} 
I_{N_{\theta}} P_{\deltap}f( n) &= \int_{N_{\theta}} \int_{\bar{P}_{\theta}} \deltap^{-1}(p) f(pn) \, dp\,dn \\
&= \int_{\bar{N}_{\theta}}\int_{M_{\theta}} \int_{N_{\theta}} \deltap^{-1}(m) f(\bar{n}mn ) dn\,dm\,d\bar{n} \\
&=\gamma(G|M_{\theta}) \int_{\bar{N}_{\theta}}\int_{M_{\theta}}\int_K  \deltap^{-1}(m)f(\bar{n}mk)\,dk\,dm\,dn \\
&=\gamma(G|M_{\theta}) \int_K P_{\deltap}(f)(k) dk.\end{align*} 
We may choose $f$ such that $ P_{\deltap}(f)(k) = \la u(k) , v(k) \ra_{\sigma}$ proving the lemma. 
\end{proof}

\subsection{}
Let $\Sigma(P_0, A_0)$ denote the set of all positive roots of $G$ with respect to $P_0$. The resulting roots by restricting to $A_{\theta}$ is denoted $\Sigma(P_{\theta}, A_{\theta})$. 

Define $\mathfrak{a}_{0, \mathbb{R}}:=(\hom_k(A_0, k^{\times})\otimes_{\mathbb{Z}} \mathbb{R})^*$ where $*$ denotes the real dual of the vector space. Letting $A_{\theta}$ denote the center of $M_{\theta}$ we may define in an analogous fashion, the vector space $\mathfrak{a}_{\theta,\mathbb{R}}.$  Then there is a canonical decomposition $$\mathfrak{a}_{0,\mathbb{R}}=\mathfrak{a}_{\theta, \mathbb{R}} \oplus \mathfrak{a}^{\theta}_{\mathbb{R}}$$ realizing $\mathfrak{a}_{\theta, \mathbb{R}}$ as a subspace of $\mathfrak{a}_{0,\mathbb{R}}$. The complexified vector spaces are denoted $\mathfrak{a}_{\theta,\mathbb{C}}$.

Let $$\mathfrak{a}_{0, \mathbb{R}}^{+}:=\{H \in \mathfrak{a}_{0,\mathbb{R}}\mid \forall \alpha \in \Sigma(P_0,A_0), \;\langle \alpha, H \rangle > 0 \}.$$  We call this the (open) positive chamber of $\mathfrak{a}_{0,\mathbb{R}}$ with respect to $P_{0}$. Similarly we define $$\mathfrak{a}^+_{\theta, \mathbb{R}}:=\{H\in \mathfrak{a}_{\theta, \mathbb{R}} \mid \forall \alpha \in \Sigma(P_{\theta}, A_{\theta}),\; \langle\alpha, H\rangle> 0\}.$$

Define the \lq $\log$\rq\ map $H_{M_{\theta}}: M_{\theta}\rightarrow \mathfrak{a}_{\theta,\mathbb{R}}$  where  $H_{M_{\theta}}(m)$ is the element of $\mathfrak{a}_{\theta, \mathbb{R}}$ such that  $$q^{-\langle \nu, H_{M_{\theta}}(m)\rangle}= |\nu(m)|$$ for all $\nu\in \hom_k(M_{\theta}, k^{\times})$. 
Let $M_{\theta}^+:=H^{-1}_{M_{\theta}}(\mathfrak{a}_{\theta,\mathbb{R}}^+)$ and $A_{\theta}^+:=M_{\theta}^+ \cap A_{\theta}$. 

We describe a partition of $A^+_{0}$. For any $0<\epsilon \leq 1$,  define \begin{align*}
A_{0}^+(\theta,\epsilon)=\{ a\in  A_0^+ \mid\quad &|\alpha(a)|\leq \epsilon\;\forall \alpha \in \Delta - \theta\\
\;\;\quad\epsilon <&|\alpha(a)|\leq 1\;\forall \alpha \in \theta\}.\end{align*}
Then we have a disjoint union $$A_0^+=\bigcup_{\substack{\theta}} A_{0}^+(\theta,\epsilon).$$  

\subsection{}
Let $\hom(G, \mathbb{C}^{\times})$ denote the group of continuous homomorphisms from $G$ to $\mathbb{C}^{\times}$. If $\chi \in \hom_k(G, k^{\times})$ then $|\chi|_k$ is defined by $|\chi|_k(g)=|\chi(g)|_k$. Let $G^1:= \bigcap_{\chi} \mathrm{Ker} |\chi|_k$ and let $X_{ur}(G):=\hom(G/G^1, \mathbb{C}^{\times})$. We call this set of characters the set of \emph{unramified} characters of $G$. These definitions apply to any Levi subgroup $M_{\theta}$ of any standard proper parabolic subgroups of $G$ by replacing $G$ with $M_{\theta}$. 

There is a surjection from $\mathfrak{a}^*_{\theta,\mathbb{C}}$ onto $X_{ur}(M_{\theta})$ defined by $\chi \otimes s \mapsto (g \mapsto |\chi(g)|^s)$. If $\nu \in \mathfrak{a}_{\theta, \mathbb{C}}^*$, the corresponding character in $X_{ur}(M_{\theta})$ is denoted $\chi_{\nu}$. The kernel of this map is a lattice of the form $(2\pi \sqrt{-1} / \log q)R$ where $R$ is a lattice of $\hom_k(M_{\theta}, k^{\times})\otimes_{\mathbb{Z}}\mathbb{Q}$. This endows $X_{ur}(M_{\theta})$ with the structure of a complex algebraic variety isomorphic to $(\mathbb{C}^{\times})^d$ where $d$ is the dimension of $\mathfrak{a}_{\theta, \mathbb{R}}$. 

Given $\chi \in X_{ur}(M_{\theta})$, suppose that $\lambda \in \mathfrak{a}^*_{\theta, \mathbb{C}}$  projects to $\chi$. We define $\Re \chi=\Re \lambda$. This is well defined as $\Re \lambda$ is independent of the choice of $\lambda$. If $\chi \in \hom(A_{\theta}, \mathbb{C}^{\times})$, then the character $|\chi|$ extends uniquely to an element of $X_{ur}(M_{\theta})$ taking positive real values. Set $\Re\chi = \Re|\chi|$. 

Consider the group $\mathrm{Im} X(A_{\theta})$ consisting of characters of $\chi\in \hom(A_{\theta}, \mathbb{C}^{\times})$ satisfying $\Re\chi =0$. Then there is a surjection with finite kernel of $\mathrm{Im} X_{ur}(M_{\theta})$ onto $\mathrm{Im}X(A_{\theta})$. The latter group is compact and we give it the Haar measure normalized so that its total mass is 1. We pull back the measure on $\mathrm{Im}X(A_{\theta})$ to a measure on $\mathrm{Im}X_{ur}(M_{\theta})$ denoting this as $d\nu$.

\newcommand{\bart}{\bar{\theta}}

\subsection{}
Let $\Xi(g)$ the \emph{zonal spherical function} of $G$ and $\Xi^M$ for the same for each standard Levi subgroups $M$ of $G$.  

Now consider a $k$-rational representation of $G$ into some $GL(n)$ with compact kernel. By pulling back to $G$, define $|| g||:=\sup\{|g_{ij}|, |(g^{-1})_{ij}|\}$ and $\sigma(g):=\log_q ||g||$ and $\sigma_*(g):=\inf_{z\in Z_G}\sigma(gz)$ where $Z_G$ denotes the center of $G$. By construction, $\sigma_*$ is subadditive and identically zero on the maximal compact subgroup $K$. 

Let $C(H_j\setminus G / H_j; \chi)$ denote the space of $H_j$ bi-invariant smooth complex valued functions on $G$ with $Z_G$ acting by a unitary character $\chi$.
For any $f \in C(H_j\setminus G/H_j)$ define $$q_{1,r}(f):= \sup_{g \in G} f(g) \Xi^{-1}(g)(1+ \sigma_*(g))^{r}$$ where $ r > 0$ .
 Define \begin{align*}\mathcal{C}^*(G;\chi):=
\bigcup_j\{&f\in C(H_j\setminus G/ H_j; \chi) \mid  q_{1,n}(f) < \infty \;\;\forall n\in \mathbb{N}\}.\end{align*}
Then $\{q_{1,n}\}_{n\in \mathbb{N}}$ is a family of continuous seminorms making $\mathcal{C}^*(G)=\mathcal{C}^*(G; \chi)$ into a Frechet space. 
It is well known that $\mathcal{C}^*(G;\chi)$ is dense (with respect to the $L^2$ norm) in $L^2(G;\chi)$. 

Let $\psi$ be a nondegenerate character of $N_0$. Define $C(N_0\setminus G/H_j; \psi)$ to be the space of right $H_j$ invariant complex valued functions of $G$ such that $f(ng)=\psi(n)f(g)$ for all $n \in N_0$ and $g\in G$. We agree that the center $Z_G$ of $G$ acts by character $\chi$ and suppress this from subsequent notation. 

Consider any $f \in C(N_{\emptyset}\setminus G/ H_k)$. By the Iwasawa decomposition $G=N_{\emptyset}M_{\emptyset}K$ write $g=n_{\emptyset}(g)m_{\emptyset}(g)k(g)$ with $m(g) \in M_{\emptyset}$, $n_{\emptyset}(g)\in N_{\emptyset}$ and $k(g) \in K$. Although $m_{\emptyset}(g)$ is only defined up to $M_{\emptyset} \cap K$, $\delta_{P_0}^{-\frac{1}{2}}(m)=1$ for any $m \in M_{\emptyset} \cap K$ implying that for any $r>0$, $$q_{2,r}(f):=\sup_{g \in G} f(g) \delta_{P_0}^{-\frac{1}{2}}(m_{\emptyset}(g))(1+\sigma_*(g))^r$$ is well defined.

Now define the space
\begin{align*}\mathcal{C}^*(N_{\emptyset}\setminus G;\psi)&:=
\bigcup_j\{f\in C(N_{\emptyset}\setminus G/ H_j;\psi) \mid q_{2,n}(f) < \infty \;\;\forall n \in \mathbb{N}\}.\end{align*}
As before the set $\{q_{2,n}\}_{n\in\mathbb{N}}$ forms a family of continuous seminorms on $\mathcal{C}^*(N_{\emptyset}\setminus G;\psi)$ making it a Frechet space. It is dense in $L^2(N_{\emptyset}\setminus G;\psi)$. 

Let $P=P_{\theta}$ be a proper standard parabolic subgroup of $G$ with $M_{\theta}$  its Levi subgroup and writing $M_{\theta} \cap N_{\emptyset}$ as $N_*$.   Write $\sigma_*^{M_{\theta}}(m):=\inf_{a\in A_{\theta}}\sigma(ma)$. This generalizes $\sigma_*$ defined earlier. To be precise, $\sigma_*=\sigma_*^G$. 

We define \begin{align*}q^{M_{\theta}}_{1,r}(f)&:= \sup_{m \in {M_{\theta}}} f(m) \Xi^{M_{\theta}}(m)^{-1}(1+\sigma^{M_{\theta}}_*(m))^r, \\  
q^{M_{\theta}}_{2,r}(f)&:=\sup_{m \in {M_{\theta}}}f(m)\delta_{P_0 \cap {M_{\theta}}}^{-\frac{1}{2}}(m_{\emptyset}(m))(1+\sigma^{M_{\theta}}_*(m))^r.\\
\intertext{and}
q^{G;M_{\theta}}_{2,r}(f)&:=\sup_{m \in {M_{\theta}}}f(m)\delta_{P_0 \cap {M_{\theta}}}^{-\frac{1}{2}}(m_{\emptyset}(m))(1+\sigma^G_*(m))^r.
\end{align*} 

Define the following Frechet spaces
\begin{align*}
\mathcal{C}^*({M_{\theta}}):=\bigcup_j\{ &f\in C(H_j\cap {M_{\theta}}\setminus {M_{\theta}}/{M_{\theta}}\cap H_j)\mid q_{1,n}^{M_{\theta}}(f) < \infty \;\; \forall n \in \mathbb{N}\}, \\
\mathcal{C}^*(N_*\setminus {M_{\theta}} ;\psi):=\bigcup_j\{&f\in C(N_*\setminus {M_{\theta}}/ {M_{\theta}}\cap H_j;\psi^{\theta}) \mid q^{M_{\theta}}_{2,n}(f) <\infty \;\;\forall n \in \mathbb{N}\},\\ 
\intertext{and}
\mathcal{C}^*(G;N_*\setminus {M_{\theta}} ;\psi):=\bigcup_j\{&f\in C(N_*\setminus {M_{\theta}}/ {M_{\theta}}\cap H_j;\psi^{\theta}) \mid q^{G;M_{\theta}}_{2,n}(f) <\infty \;\;\forall n \in \mathbb{N}\}
\end{align*} where $\psi^{\theta} =: \mathrm{Res}^{N_{\emptyset}}_{N_*}(\psi)$.  
Then we have an inclusion $$\mathcal{C}^*(G;N_*\setminus {M_{\theta}}; \psi^{\theta}) \subset \mathcal{C}^*(N_*\setminus {M_{\theta}}; \psi^{\theta}).$$
We will refer to these spaces generically as \emph{Schwartz spaces}. 

\section{Whittaker functions and the Harish-Chandra transform}

\subsection{}
Let $(\pi,V)$ be a smooth finitely generated representation of $G$ and $V'$ denote the algebraic complex linear dual  of $V$. Define $(\pi'(g)\lambda)(v)=\lambda(\pi(g^{-1})v)$ and let $V^{\vee}$ denote the smooth points of ${\pi'}$. Let $\pi^{\vee}$ be the representation obtained by restricting $\pi'$ to $V^{\vee}$. Then $(\pi^{\vee}, V^{\vee})$ is a smooth representation of $G$. 

Now consider a fixed $\lambda \in V'$ and write $$W(v,\lambda)(g):=\lambda(\pi(g)v)$$ with $v\in V$. Since $V$ is smooth, $W(v, \lambda)$ is $K_j$ invariant on the right for some open compact subgroup. 
Note that if $\check{v}\in V^{\vee}$, then $W(v,\check{v})(g)$ is  a matrix coefficient of $V$. 

Fix \emph{any} (not necessarily nondegenerate) character  $\psi$ of $N_0$. Fix a standard parabolic subgroup $P_{\theta}=P=MN$ of $G$ corresponding to $\theta \subset \Delta$ ($P$ possibly equals to $G$) for the discussion throughout this subsection. Let $N_*$ and $\psi^{\theta}$ be the unipotent subgroup and character defined as in the previous section. Then define $$V(N_*,\psi^{\theta}):=\mathrm{span}\{\pi(n)v-\psi^{\theta}(n)v \mid v\in V, \;\; n\in N_* \}$$  and $$V(N):=\mathrm{span}\{\pi(n)v-v \mid v\in V ,\;\; n\in N \}.$$ 

Define vector spaces $$V_{N_*, \psi^{\theta}}=V/V(N_*, \psi^{\theta})$$ and $$V_{N} = V/V(N).$$ Then $$r^{G}_{P}(V)=V_{N}\otimes \delta_{P}^{-1/2}$$ is the Jacquet restriction functor sending smooth, finitely generated admissible representations of $G$ to smooth, finitely generated admissible representations of $M$.    

Now we require $\psi$ to be a nondegenerate character. Define $$Wh_{\psi}(V):= \hom_{N_0}(\pi, \mathbb{C}_{\psi}).$$

A representation $\pi$ is said to be $\psi$-generic if $Wh_{\psi}(V)$ is nontrivial. Now assume that this  is the case for $\pi$ and consider any nonzero $ \lambda \in Wh_{\psi}(V)$ and nonzero $v \in V$.  Then $W(v, \lambda)(g)$ is not identically zero as a function on $G$ and satisfies $$W(v,\lambda)(ng)= \psi(n)W(v,\lambda)(g)$$ for any $g \in G$ and $n \in N_0$. 

We say that $W(v,\lambda)$ is a \emph{Whittaker function} and  $W(v,\lambda)(1)=\lambda(v)$  a \emph{Whittaker functional}. All these can be generalized to Levi subgroups of $G$ as well, by substituting $G$ for $M$, $N_0$ for $N_*$ and $\psi$ for $\psi^{\theta}$. 

Let $\Phi_{\theta}$ be the canonical map from $(V_{N_{\theta}})_{N_*, \psi^{\theta}}$ to $V_{N_{0},\psi}$ introduced in \cite{casselman} (see also \cite{casselnotes}) . If $v\in V$, we write $\tilde{v}$ for its image in $(V_{N_{\theta}})_{N_*, \psi^{\theta}}$. 
Recall the following lemma of Casselman. (See \cite[Proposition 6.3 and 6.4]{casselman}).

\begin{plemma}\label{partition of whit func} Fix a standard parabolic subgroup $P_{\theta}$ of $G$. Let $\lambda\in Wh_{\psi}(V)$ and $v\in V$. Then there exists $\epsilon >0$ such that $$W(v, \lambda)(a)=W(\tilde{v}, \lambda \circ \Phi_{\theta})(a)$$ for any $a\in A_0$ satisfying $|\alpha(a)|< \epsilon  \;\; \forall \alpha \in \Delta -\theta.$
\end{plemma}
Notice that when $\theta=\emptyset$, this is Proposition 6.3 in \cite{casselman}. 


\subsection{}
For a fixed but arbitrary $\theta \subset \Delta$, consider an $A_{\theta}$-finite complex valued smooth function $f$ on $A_{\theta}$. Then it is well known  that $$f(a)=\sum^{}_{\substack{\nu}} \nu(a)P(H_{M_{\theta}}(a))$$ where $P(x)$ is a polynomial on the Lie algebra $\mathfrak{a}_{\theta, \mathbb{R}}$ of $A_{\theta}$ and $\nu$ a smooth  character of $A_{\theta}$. The characters occurring in the decomposition above are known as the \emph{exponents} of $f$. Let $E(A_{\theta}, f)$  denote the set of exponents of $f$. 
 
\begin{plemma}\cite[Proposition 4.4.4]{casselnotes}\label{summability in the positive cone}
Let $\theta \subset \Delta, \epsilon \in (0,1]$ and $p>0$. Let $f:A_{0}^+(\theta,\epsilon)\rightarrow \mathbb{C}$ be a complex valued function such that \begin{enumerate}
\item $f$ is the restriction to $A_{0}^+(\theta,\epsilon)$ of an $A_{\theta}$-finite function; 
\item the center of $G$, $Z_G$ acts by a unitary character on $f$ and 
\item $f$ is invariant under right translation by some open subgroup $A_{K_i}$ of $ A_{0}^1.$\end{enumerate}
Then $|f|^p$ is integrable on $A_{0}^+(\theta,\epsilon)/A_{K_i}Z_G$ if and only if $|\chi(a)| <1$ for all $a \in A^+_{\theta}$ and each $\chi \in E(A_{\theta}, f)$. 
\end{plemma}

It is well known that if a representation of $G$, $(\pi,V)$ is smooth and admissible, then so is $V_N$. Thus the center of $M$, $A$, acts locally finitely on $V_N$ so that  $$V_N=\bigoplus_{\nu} (V_N)_{\nu}$$ as generalized eigenspaces where $\nu$ are smooth  characters of $A$. The set of all characters appearing in this decomposition are called $\emph{exponents}$ of $V$ with respect to $P$ denoted $E(P, V)$. If $V$ is finitely generated, this set is finite. 

An irreducible smooth representation is said to be a \emph{discrete series} (resp. \emph{tempered}) representation of $G$ if the center acts by a unitary character and its matrix coefficients (modulo the center) are in $L^2(G)$ (resp. $L^{2+\epsilon}(G)$ for any $\epsilon > 0$). The following result is well known (see \cite[Proposition III.1.1 and Proposition III.2.2]{walsp}). 

\begin{prop}\label{exponents of discrete series}
Suppose $(\pi,V)$ is a discrete series representation (resp. tempered representation) of $G$. Then it is necessary and sufficient that for every standard parabolic subgroup $P$ and every $\nu\in E(P,V)$, $|\delta^{-\frac{1}{2}}_P(a)\nu(a)| < 1$ (resp. $ |\delta^{-\frac{1}{2}}_P(a)\nu(a)| \leq 1$)for all $a\in A^+$.  
\end{prop}
\begin{flushright}
$\qed$ 
\end{flushright}

\begin{pthm}\label{suff}
Let $\psi$ be an nondegenerate additive unitary character of $N_0$ and let $(\pi,V)$ be an irreducible $\psi$-generic discrete series representation, then $V$ embeds into  $L^2(N_0\setminus G;\psi)$. 
\end{pthm}

\begin{proof}
Any Whittaker function is of the form $W(v,\lambda)(g)$ for $\lambda\in (V_{N_0,\psi})'$. Now choose $K_i$ as the largest open compact subgroup of $K$ such that $$\int_{K_i}\pi(k)v\; dk=v'$$ does not vanish. Let $[K / K_i]$ denote a (finite) set of coset representatives of $K / K_i$. Then for each $k_j \in [K / K_i]$, set $v_j'=k_jv'.$  

Then,  \begin{align*}
\int_{N_0\setminus G/Z_G} |W(g)|^2 d\bar{g}&=\int_{A_0/A_0^1Z_G}\int_{K} \delta_0^{-1}(a)|W(ak)|^2 da dk \\
&=\int_{A_0/A_0^1Z_G}\sum_{k_j\in [K/ K_i]}\int_{K_i} \delta_0^{-1}(a)|W(ak_jk_i)|^2 da dk_i \\
&=\int_{A_0/A_0^1Z_G}\sum_{j} \delta_0^{-1}(a)|W(v_j',\lambda)(a)|^2 da. \\
\end{align*}

Thus, the finiteness of that integral depends upon the square integrability of $W(v_j',\lambda)(a)$ on $A_0/A_0^1Z_G$. Since $\psi$ is nondegenerate, we observe that $W(v_j',\lambda)(a)$ is supported inside a translate of  $A_0^+$ (\cite[Proposition 6.1]{casselman}). By replacing $v_j'$ with $\pi(a_j)v_j'$ for some suitable $a_j \in A_0$, we may even assume that this support is contained in $ A_0^+$.

Choose $\epsilon=\min\{\epsilon_{\theta}\}_{\theta}$ where $\epsilon_{\theta}$ is obtained by applying Lemma \ref{partition of whit func} to $P_{\theta}$ (with $\theta\subset \Delta$). Recall that $A^+_0$ is partitioned into $\bigcup_{\theta} A_{0}^+(\theta,\epsilon)$. We let $f(a):=\delta_0^{-1/2}(a).W(v',\lambda)(a)$ and restrict $f(a)$ to each partition so that we may then apply Lemma \ref{summability in the positive cone} with $p=2$.

We must check all three conditions of this lemma. Condition (2) is clear and (3) is satisfied by $A_{K_i}=A_0^1$. 
As $\pi$ is admissible, so is $V_N$ and thus $A_{\theta}$ acts locally finitely on $(V_{N_{\theta}})_{N_*, \psi^{\theta}}$. Therefore $W(\tilde{v}^{\prime}, \lambda \circ \Phi_{\theta})(a)$ satisfies condition (1). 

Clearly the exponents of  $W(\tilde{v}^{\prime},\lambda\circ \Phi_{\theta})(a)$ are in $E(P_{\theta},V).$ By Proposition \ref{exponents of discrete series}, this will imply the square integrability of $W(v', \lambda)$ in $A_0^+$. 
This proves the theorem. 
\end{proof}

We introduce the notion of moderate growth for functions on $N_0 \setminus G$. Define $\mathscr{A}(N_0\setminus G ;\psi)$ to be the union over all open compact subgroups $K_i$ of $G$ of subspaces of functions in $C(N_0\setminus G/K_i; \psi)$ defined by the following growth condition: There exists a constant $C > 0$ and  $r>0$ such that $$|f(g)| \leq C \delta^{\frac{1}{2}}_0(m_0(g))(1+\sigma_*(g))^r.$$ 

\begin{prop} \label{tempered}
Assume that $(\pi,V)$ is a tempered and $\psi$-generic irreducible representation of $G$. Then for any $v \in V$ and $\lambda \in Wh_{\psi}(V)$, $W(v, \lambda) \in \mathscr{A}(N_0 \setminus G; \psi)$.  
\end{prop}
\begin{proof}
We know that if $\pi$ is tempered then for each standard parabolic subgroup $P_{\theta}$ and exponent $\nu \in E( \theta,V)$,  $|\delta_P^{-\frac{1}{2}}(a)\nu(a)|\leq 1$ by \cite[Proposition III.2.2]{walsp}. If this inequality was strict for all parabolic subgroups and all exponents, then $\pi$ is a discrete series and hence by (the proof of) Theorem \ref{suff}, $W(v, \lambda)$ is in $\mathscr{A}(N_0 \setminus G; \psi)$. Thus we assume that for some parabolic subgroup, $P_{\theta}$, there exists an exponent such that $|\delta_P^{-\frac{1}{2}}(a)\nu(a)|= 1.$ In this case, using Lemma \ref{partition of whit func} and arguing as in Theorem \ref{location of image} we obtain the result. 
\end{proof}


\subsection{}
If $P$ is a parabolic subgroup of $G$, then $\bar{P}$ will denote its opposite parabolic subgroup. Given $f\in \mathcal{C}^*(N_0\setminus G;\psi)$ define the \emph{Harish-Chandra} transform \begin{align} f^{P}(m)=\delta_P^{1/2}(m)\int_{\bar{N}}f(\bar{n}m)d\bar{n}=\delta^{1/2}_{\bar{P}}(m)\int_{\bar{N}}f(m\bar{n})d\bar{n} \label{hcdef}.\end{align}

\begin{lemma}\label{switching order}
Let $P_{\theta}$ be a standard parabolic subgroup and $N_{\theta}$ its unipotent subgroup. For any $f \in \mathcal{C}^*(G)$, there exists a continuous seminorm $q(f)$ such that $$ \int_{N_0 \times \bar{N}_{\theta}} | f(n_0 \bar{n})| \, dn_0 d\bar{n} < q(f).$$
\end{lemma}
\begin{proof}
Using the Iwasawa decomposition, write $\bar{n}= n_0(\bar{n}) a_0(\bar{n}) k(\bar{n})$. Then the integral can be written as $$ \int_{N_0 \times \bar{N}_{\theta}} | f(n_0 a_0(\bar{n}) k(\bar{n})|\, dn_0 d\bar{n}. $$ By \cite[Proposition II.4.5]{walsp}, we know that for any integer $d > 0$, one can find and integer $r >0$ (and seminorm $q_{1,r}(f)$) such that this integral is majorized by $$ \int_{\bar{N}_{\theta}} \delta_0(a_0(\bar{n}))^{\frac{1}{2}} (1+ \sigma_*( a_0(\bar{n})))^{-d}\, d\bar{n}.$$  But we know by \cite[Lemme II.4.2]{walsp}  that this integral converges for sufficiently large $d$. 
\end{proof}

\begin{lprop}\label{corrected}
The integral in \eqref{hcdef} converges absolutely and uniformly over compact sets in $M$, $f^P\in \mathcal{C}^*(N_0\cap M\setminus M;\psi)$ and $f\mapsto f^P$ is continuous in the topology induced by seminorms on $\mathcal{C}^*(N_0\setminus G;\psi)$ and $\mathcal{C}^*(G;N_0\cap M\setminus M;\psi)$. \end{lprop}

\begin{proof}
Write $\bar{n}m = n(\bar{n}m) a(\bar{n}m) k(\bar{n}m)$ using the Iwasawa decomposition. 
If we write $m\in M$ as $m = n_1 m_1 k_1$ where $n_1 \in N_0 \cap M$, $m_1 \in A_0 \cap M$ and $k_1 \in K \cap M$, then $$\bar{n}m = \bar{n}n_1 m_1 k_1 = n_1m_1.( (n_1m_1)^{-1} \bar{n} (n_1m_1)) k_1.$$ This implies that $$a(\bar{n}m)= m_1 a(( n_1m_1)^{-1} \bar{n} (n_1m_1))$$ and $$n(\bar{n}m)=n_1m_1^{}n((n_1m_1)^{-1}\bar{n}(n_1m_1))m_1^{-1}.$$ 

Let $\bar{\mathbf{n}}=(n_1m_1)^{-1}\bar{n}(n_1m_1)$ so that $$n(\bar{n}m)a(\bar{n}m)=n_1m_1n(\bar{\mathbf{n}})a(\bar{\mathbf{n}}).$$ Note that $n_1m_1 \in M$ and $\bar{\mathbf{n}} \in \bar{N}$. 

By hypothesis, we fix a compact set $S\subset M$ and prove the convergence of $f^P(m)$ as $m$ varies over $S$. Let $C = \max(1+\sigma_*(n_1m_1))^d.$ Then we claim that $$ 1 \leq \left(\frac{1+\sigma_*(n_1m_1n(\bar{\mathbf{n}})a(\bar{\mathbf{n}}))}{1+\sigma_*(n(\bar{\mathbf{n}})a(\bar{\mathbf{n}}))}\right)^d \leq C$$ for all $\bar{n} \in \bar{N}$ with $\sigma_*(\bar{\mathbf{n}})$ sufficiently large. 

Indeed, the right hand side of the inequality follows easily from the subadditivity of $\sigma_*$. To prove the left hand side of the inequality we first note that  $$\sigma_*(n_1m_1n(\bar{\mathbf{n}})a(\bar{\mathbf{n}}))=\sigma_*(n_1m_1n(\bar{\mathbf{n}})a(\bar{\mathbf{n}})k(\bar{\mathbf{n}}))=\sigma_*(n_1m_1\bar{\mathbf{n}})$$ as $\sigma_*$ is $K$ invariant. Now by a lemma of Waldspurger (\cite[Lemme II.3.1]{walsp}), $1+\sigma_*(n_1m_1\mathbf{\bar{n}})$ is bounded below by $1+ \sup\{\sigma_*(n_1m_1), \sigma_*(\bar{\mathbf{n}})\}.$ Since $\sigma_*(n_1m_1)$ is bounded, the claim follows immediately. 

With this, we may choose a suitable positive constant $A$, so that $$C \leq A\left(\frac{1+\sigma_*(n_1m_1n(\bar{\mathbf{n}})a(\bar{\mathbf{n}}))}{1+\sigma_*(n(\bar{\mathbf{n}})a(\bar{\mathbf{n}}))}\right)^d$$ implying that $$(1+\sigma_*(n_1m_1))^d(1+\sigma_*(n_1m_1n(\bar{\mathbf{n}})a(\bar{\mathbf{n}})))^{-d} \leq A (1+\sigma_*(n(\bar{\mathbf{n}})a(\bar{\mathbf{n}})))^{-d}.$$

A function $f \in \mathcal{C}^*(N_0\setminus G ; \psi)$ satisfies $$|f(\bar{n}m)| < \delta^{\frac{1}{2}}_0(a(\bar{n}m)) (1+ \sigma_*(n(\bar{n}m)a(\bar{n}m)))^{-d}q_{2,d}(f).$$ Then from the discussion above, we obtain
 \begin{align*} \delta^{-\frac{1}{2}}_0(m_1) (1+ \sigma_*(n_1m_1))^{d}|f(\bar{n}m)| &< \delta_{0}^{\frac{1}{2}}(a(\bar{\mathbf{n}}))(1+ \sigma_*(n(\mathbf{\bar{n}}) a(\mathbf{\bar{n}})))^{-d}q_{2,d}(f)\\
  &= \delta_{0}^{\frac{1}{2}}(a(\bar{\mathbf{n}}))(1+ \sigma_*(\mathbf{\bar{n}}))^{-d} q_{2,d}(f). 
  \end{align*}
For $m \in S$ \begin{align}  \delta^{-\frac{1}{2}}_0(m_1) (1+ \sigma_*(n_1m_1))^{d}\int_{\bar{N}_{\theta}} |f(\bar{n}m)| \,d\bar{n}\label{reff2}\end{align} is majorized by a constant multiple of
 \begin{align} \int_{\bar{N}_{\theta}} \delta_0^{\frac{1}{2}}(a(\mathbf{\bar{n}})) (1+ \sigma_*(\mathbf{\bar{n}}))^{-d}\, d\mathbf{\bar{n}}= \delta_P^{-1}(m_1) \int_{\bar{N}_{\theta}} \delta_0^{\frac{1}{2}}(a(\bar{n})) (1+ \sigma_*(\bar{n}))^{-d} \,d\bar{n}.\label{reff3}\end{align}
  By \cite[Lemme II.4.2]{walsp}, the right hand side is finite for some sufficiently large $d$. This proves the first  assertion of the lemma. 
  From \eqref{reff2} and \eqref{reff3}, we have also that $$ \delta_P^{\frac{1}{2}}(m)|f^P(m)|< c.\delta^{\frac{1}{2}}_{M\cap P_0}(m_1)(1+  \sigma_*(n_1m_1))^{-d}<c'.\delta^{\frac{1}{2}}_{M\cap P_0}(m_1)(1+  \sigma_*(m))^{-d} $$ where we recall that $m=n_1m_1k_1$ proving the second assertion. Now the continuity of the Harish-Chandra transform is immediate from this proof.  
\end{proof}

Define \begin{align*}\,^{\circ}\mathcal{C}^*(N_0\setminus G;\psi)=\{ &f\in\mathcal{C}^*(N_0\setminus G;\psi)| (R(k)f)^P\equiv 0 \\ &\;\;for\;all \;k\in K\; \;and\;for\;\,all\;\;parabolic\;\;subgroups\;\;P=P_{\theta},\; \;\theta \subsetneq \Delta\}.\end{align*} We refer to this space as the space of \emph{discrete} functions. By the Iwasawa decomposition this space is stable under right translation by the full group $G$.

Let $f\in \mathcal{C}^*(N_0\setminus G;\psi)$ and $V_f$ be the space spanned by right translates of $G$ and $V_{f^P}$ the space generated by right translates of $M$ acting on $f^P$ in $\mathcal{C}^*(G;N_0\cap M\setminus M;\psi)$. 
\newcommand{\n}{\bar{n}}
\begin{lemma}\label{quotient of jacquet}
The representation space $V_{f^P}$ is a quotient of $r^G_{\bar{P}}(V_f)$. 
\end{lemma}
\begin{proof}It is easy to check from definition that $(R(m)f)^{P}= \delta_{\bar{P}}^{\frac{1}{2}}(m)R(m)f^{P}$. 

Next, let $\bar{N}$ be the unipotent subgroup opposite of $P$ and consider the space $V_f(\bar{N})$. If $h \in V_f(\bar{N})$, then $$ \int_{\bar{N}_0} \; R(n)h=0$$ for some open compact subgroup $\bar{N}_0$ of $\bar{N}$. In particular, $V_f(\bar{N})$ is contained in the kernel of $f \mapsto f^P$. The lemma follows from these two observations. 
\end{proof}
\begin{lthm}\label{location of image}
Let $(\pi,V)$ be a smooth irreducible representation of $G$ which embeds into $L^2(N_0\setminus G; \psi)$. Then the image of $V$ under this embedding is in $\;^{\circ}\mathcal{C}^*(N_0\setminus G;\psi)$. 
\end{lthm}
\begin{proof}
Let $T$ be the embedding map and define $W(v,T)(g)=T(\pi(g)v)(1)$. Then $W(v,T)$ is a Whittaker function (with $(Tv)(1)$ considered a Whittaker functional). We note that $W(v,T)$ is smooth while $\psi$ is nondegenerate and therefore supported in a translate of $ A_0^+$. 
By Lemma \ref{partition of whit func}, the restriction of $W(v,T)$ on each partition $A_{0}^+(\theta,\epsilon)$ coincides with some $A_{\theta}$-finite function. Thus conditions (1) to (3) of Lemma \ref{summability in the positive cone} are satisfied 
and we may conclude that the exponents (of the $A_{\theta}$-finite function) must satisfy $|\delta^{-\frac{1}{2}}_0(a)\nu(a)|<1$ for every $a\in A_{0}^+(\theta,\epsilon)$ as $W(v,T)$ is also square integrable. 

Since 
 $$W(v,T)(a)=\tilde{W}(\tilde{v},T\circ \Phi_{\theta})(a)=\sum_{\nu}\nu(a)P(H_{M_{\theta}}(a))$$ 
and we know that $|P(H_{M_{\theta}}(a))| \in O((1+\sigma_*(a))^n)$ for some integer $n$  (see \cite[pg. 242]{walsp}), there exists a constant $C>0$ such that 
 $$| \delta_0^{-\frac{1}{2}}(a)W(v,T)(a)| < C(1+ \sigma_*(a))^{-r})$$ for all positive integers $r$. Applying this same argument to each partition will imply that $W(v,T)$ is in $\mathcal{C}^*(N_0\setminus G;\psi)$ by definition of the Schwartz space. 

 Now Jacquet restriction preserves admissibility. Thus by Lemma \ref{quotient of jacquet}, $W(v,T)^P$ is an $A$-finite function. However, we know that $W(v,T)^P \in \mathcal{C}^*(G;N_0\cap M\setminus M;\psi).$ This certainly contradicts the decay properties of $W(v,T)^P$ unless $W(v,T)^P$ is identically zero or $P=G$. This proves the theorem. 
\end{proof}


\subsection{}
We end this section by giving an application of this theory to automorphic forms. In this, we will extend a result of Savin, Khare and Larsen \cite{sav}. By an integrable discrete series representation of $G$, we mean a discrete series representation which has matrix coefficients  in $L^1(G) \cap L^2(G)$. 
\begin{lthm} \label{main2}
Let $G$ be an almost simple, quasi-split algebraic group defined over $\mathbb{Q}$. Fix two finite and disjoint sets of places: $D$ contains $\infty$ and perhaps nothing else, and $S$ is a non-empty set of primes such that $G$ is unramified at all primes $p$ not contained in $D\cup S$. Let $\psi$ be a nondegenerate character of $N_0(\mathbb{A})$ trivial of $N_0(\mathbb{Q})$. Suppose that $\pi'_{\infty}$ is a $\psi$-generic integrable discrete series of $G(\mathbb{R})$, \textbf{and} \textbf{for every} $q$ in $D$, $\pi'_q$ is a $\psi$-generic \textbf{square integrable representation}, then there exists a global $\psi$-generic cuspidal representation $\pi$ such that $\pi_{\infty}\cong \pi'_{\infty}$, $\pi_q\cong \pi'_q$ for every $q$ in $D$ and $\pi_p$ is unramified for every $p$ outside $D\cup S$. 
\end{lthm}
\begin{proof}
We observe that the crux of the matter in the proof of Theorem 4.5 in paper \cite{sav} is the fact that $W(v,\lambda)$ satisfy certain bounds, namely there exists a $C$ such that 
$$|W(v,\lambda)(g)|\leq  C \delta_0^{1/2}(a(g))(1+\sigma_*(g))^{-r}\;\;\forall r \geq 0.$$
or in other words, the Whittaker functions are located in the Schwartz space of $N_0\setminus G$. However, this follows immediately from Theorems \ref{suff} and \ref{location of image}. 
\end{proof}

\newcommand{\ind}{I_{\bar{P_{\theta}}}^{G} \;\sigma}
\newcommand{\indj}{(I_{P_{\theta}}^{G} \;\sigma)_{\psi, N_0}}
\newcommand{\sigmaj}{\sigma_{\wpsi, N_0\cap M_{\theta}}}
\newcommand{\psiinv}{\psi^{-1}}
\newcommand{\wpsi}{\;^{w}\psi}
\newcommand{\aplms}{}
\newcommand{\ds}{d\sigma_{\alpha, \mu}}

\section{The discrete spectrum of $L^2(N_0\setminus G;\psi)$}

\subsection{}\label{sec42}

Take any $f \in \mathcal{C}^*(G)$  and consider the \emph{Harish-Chandra} transform of $f$ defined by the following (absolutely convergent) integral: 
$$ f^P(m)=\delta_P^{\frac{1}{2}}(m)\int_{\bar{N}} f(\bar{n}m)\, d\bar{n}$$ where $P=MN$ is a standard parabolic subgroup of $G$.  Then it is well known that $f^P$ is in $\mathcal{C}^*(M)$.

We give a brief overview of the decomposition of $\mathcal{C}^*(G)$ into its various constituents. 
It is known that  $\mathcal{C}^*(G)=\mathcal{C}^*_{discr}(G)\bigoplus \mathcal{C}^*_{cont}(G)$ where $\mathcal{C}^*_{discr}(G)$ is the space of functions whose Harish-Chandra transform $f^P$ is identically zero for all proper standard parabolic subgroups, $P$. 

To describe the space $\mathcal{C}^*_{cont}(G)$ (see also \cite[pg. 302]{walsp} ), let $M_{\theta}$ be the Levi subgroup of a standard parabolic associated to $\theta \subset \Delta$.
Consider the smooth normalized induced representations of $G$, $(\pi_{\bar{P}_{\theta}, \sigma, \nu}, I_{\bar{P}_{\theta}}^{G} \sigma \boxtimes \chi_{\nu})$ where $\chi_{\nu} \in X_{ur}(M_{\theta})$ ($\nu \in \mathfrak{a}^*_{\theta,\mathbb{C}}$) and $\sigma$ is a fixed discrete series representation of $M_{\theta}$. When the context is clear, we write $I(\sigma)_{\tau}$ to denote the space in  $\mathrm{Ind}_{M_{\theta} \cap K}^K\, \sigma_{ \mid_{M_{\theta} \cap K}}$ where $K$ acts according to a finite dimensional representation  $\tau$. Then take $w \in I(\sigma)_{\tau}$ and define a function $w_{\nu}$ on $G$ in the following manner: 

By the Iwasawa decomposition $G=\overline{N_{\theta}}M_{\theta}K$. Writing $g=\bar{n}mk$, we define $$w_{\nu}(\bar{n}mk) = (\chi_{\nu}.\delta_{\bar{P}}^{\frac{1}{2}})(m)\sigma(m)w(k).$$ Thus each $w_{\nu}$ is an element of $I_{\bar{P}_{\theta}}^{G} \sigma \boxtimes \chi_{\nu}$ .

 \newcommand{\dsig}{d\mu_{\mathcal{O}}(\nu)}

Take any $\chi_{\nu} \in \mathrm{Im} X_{ur}(M_{\theta})$.
 If $\sigma$ is a discrete series representation of $M_{\theta}$, then so is $\sigma \boxtimes \chi_{\nu}.$ This defines an action of the compact abelian group $\im X_{\ur}(M_{\theta})$ on the set of all discrete series representations of $M_{\theta}$. This action has a finite stabilizer and transfers the induced measure on the quotient of $\im X_{\ur}(M_{\theta})$ by $\mathrm{Stab}_{\im X_{\ur}(M_{\theta})}(\sigma)$ to the orbits. An orbit under this action is denoted as $\mathcal{O}$. 
 
Recall the constant $j(\sigma\boxtimes \chi_{\nu})$ defined in \cite[pg. 285]{walsp}. Let $$\mu(\sigma, \nu):= j(\sigma \boxtimes \chi_{\nu})^{-1} \prod_{\alpha \in \,^{nd} \Sigma(P,A)} \gamma(M_{\alpha}|M)^2$$ with $M_{\alpha}$ defined as in Section \ref{sec32}. We call the normalized measure $$\dsig = \mu(\sigma,\nu)\, d\nu$$ the \emph{Plancherel measure} on $\mathcal{O}$.

 Fix an orbit $\mathcal{O}$ and let $\sigma$ be any representation in it. Let $w, v \in I(\sigma)$ and let $\langle , \rangle_{I(\sigma)}$ be the unitary Hermitian form on $\ind$. Then $\varphi_{w,v}[g](\nu):=\langle \pi_{\bar{P}_{\theta}, \sigma,\nu}  (g)v_{\nu}, w_{\nu} \rangle_{I(\sigma)}$ is considered as a function on $\mathcal{O}$. 

Take any infinitely differentiable  function $\alpha(\nu) \in C^{\infty}(\mathcal{O})$ and define the following integral transform from $C^{\infty}(\mathcal{O})$ to $\mathcal{C}^*(G)$:
 $$ f_{[\alpha, \mathcal{O}, \bar{P}_{\theta}];v,w}(g)=\int_{\mathcal{O}} \varphi_{v,w}[g](\nu) \alpha(\nu)\,\dsig.$$ This is the Harish-Chandra wave packet map.
 
By \cite[Theoreme VIII.1.1]{walsp}, any $f\in \mathcal{C}_{cont}^*(G)$  is a finite linear combination of $f_{[\alpha, \mathcal{O}, \bar{P}_{\theta}];v,w}$  over suitable choices of datum $[\alpha, \mathcal{O}, \bar{P}_{\theta}]$ with $P_{\theta} \neq G$ and $w,v \in I(\sigma)$.

\subsection{}
We write $\langle\;,\rangle$ for the sesquilinear inner product on $L^2(N_0\setminus G;\psi)$. For $f\in \mathcal{C}^*(G)$ and $\varphi \in \mathcal{C}^*(N_0 \setminus G; \psi)$, we let $$ B(f, \varphi) = \int_G f(g) \overline{\varphi(g)}\, dg.$$

\begin{lemma} \label{lemma42}The integral defining the sesquilinear pairing $B(f, \varphi)$ converges absolutely for all $f \in \mathcal{C}^*(G)$ and $\varphi \in \mathcal{C}^*(N_0 \setminus G ; \psi).$ Moreover there exists two continuous seminorms $q_1$ and $q_2$ on $\mathcal{C}^*(G)$ and $\mathcal{C}^*(N_0 \setminus G ; \psi)$  respectively such that $$ | B(f, \varphi)| \leq q_1(f) q_2(\varphi).$$
\end{lemma}
\begin{proof}
Now for any positive integer $d'$ and $d$, \begin{align} 
|B(f,\varphi)| &\leq \int_{N_0 \times A_0 \times K} \delta_{P_0}^{-1}(a)|f(nak)| |\varphi(ak)|\, dn\, da\, dk \notag \\
&\leq  \int_{N_0 \times A_0 \times K}\delta_{P_0}^{-\frac{1}{2}}(a)  q_{2,d}(R(k)\varphi)q_{1,d'}(R(k)f) \Xi(na)(1+\sigma_*(na))^{-d'} 
(1+\sigma_*(a))^{-d}.  \label{step1} \\
\intertext{By \cite[Proposition II.4.5]{walsp}, for any $r \geq 0$ one can find a positive constant and positive integer $d'$ such that \eqref{step1} is}
& \leq \int_{A_0 \times K} q_{1,r}(R(k)f) q_{2,d}(R(k)\varphi) (1+\sigma_*(a))^{-r-d}\, da\,dk\,\label{step2} . 
\intertext{Let $q_1(f)$ and $q_2(\varphi)$ be the supremum of $q_{1,r}(R(k)f)$ and $q_{2,d}(R(k)\varphi)$ over the entire maximal compact subgroup $K$. Then the expression in line \eqref{step2} }
& \leq q_1(f)q_2(\varphi)\int_{A_0} (1 + \sigma_*(a))^{-r-d} \,da \notag.
\end{align} However this is finite if $r+d$ is sufficiently large. 
\end{proof}

 We let $$f_{\psi}(g):=\int_{N_0}\psi(n)^{-1}f(ng)dn.$$ By the proof of Lemma \ref{lemma42} above this integral converges absolutely and  the map \begin{align*}\mathcal{C}^*(G) &\rightarrow \mathcal{C}^*(N_0 \setminus G ; \psi) \\
 f &\mapsto f_{\psi}
 \end{align*} is continuous and it is clear that the sesquilinear pairing $B(f,\varphi)$ is $G$ invariant and equals $\langle f_{\psi}, \varphi\rangle.$
  
 The following is the key lemma of the paper. 
\begin{lemma}\label{pairing}
With notation as in this section, 
\begin{enumerate} 
\item The map $f\mapsto f_{\psi}$ is surjective from $\mathcal{C}^*(G)$ to $\mathcal{C}^*(N_0\setminus G;\psi).$
\item If $f \in \mathcal{C}^*_{cont}(G)$, then $B( f, \varphi) =0$ for any $\varphi \in \,^{\circ}\mathcal{C}^*(N_0 \setminus G;\psi)$. 
\item
If $B( f,\varphi)=0$ for all $f\in \mathcal{C}^*(G)$, then $\varphi$ is identically zero. \end{enumerate}\end{lemma}
\begin{proof}
Assertions (1) and (2) will be proven in the next subsection. (3) follows easily from (1). 
\end{proof}

\begin{lthm}\label{main3}
The space of discrete functions $\;^{\circ}\mathcal{C}^*(N_0\setminus G;\psi)$ is a multiplicity free direct sum of discrete series representations of $G$ which are $\psi$-generic. 
\end{lthm}
\begin{proof}
We first show that $\,^{\circ}\mathcal{C}^*(N_0 \setminus G;\psi)\supset \mathcal{C}^*_{discr}(G)_{\psi}$.
Suppose that $f \in \mathcal{C}_{discr}^*(G)$. Then we claim $(f_{\psi})^P=0$ for any proper standard parabolic $P=P_{\theta}$. Indeed, by Lemma  \ref{switching order},  $$ \delta_P^{\frac{1}{2}}(m)\int_{\bar{N}_{\theta}}\left( \int_{N_0} \psiinv(n) f(\bar{n}mn)\, dn_0\right) \, d\bar{n} = \int_{N_0} \psiinv(n)\left( \delta^{\frac{1}{2}}_P(m) \int_{\bar{N}_{\theta}} (R(n)f)(\bar{n}m)\right) d\bar{n}.$$
Since $\mathcal{C}^*_{discr}(G)$ is stable under right translation, the right hand side is zero. 

We know that $C^*_{cont}(G)$ maps into the orthogonal complement of $\, ^{\circ}C^*(N_0 \setminus G; \psi)$. As the projection of $\mathcal{C}^*(G)$ to its $\psi$-coinvariants is surjective, $\,^{\circ}\mathcal{C}^*(N_0 \setminus G;\psi)=\mathcal{C}^*_{discr}(G)_{\psi}$. 

Let $(\sigma,H_{\sigma})$ be a discrete series representation of $G$. 
Since $G$ is assumed to be quasi-split, a result of Shalika assures us that $\dim\,Wh_{\psi}(H_{\sigma})\leq 1$. The multiplicity free assertion follows immediately from this and proves the theorem.

\end{proof}


\subsection{}
In this section and the next, we will complete the proof of Lemma \ref{pairing}. 

Suppose we are given $f\in \mathcal{C}^*(N_0\setminus G;\psi)^H$ where $H=H_m$ is some open compact subgroup.  Then $f$ is determined by the values it takes on the cosets of $N_0\setminus G/H$. More precisely, if $g\in \gamma g_iH$, then $$f(g)=\begin{cases} \psi(\gamma)f(g_i) &\; \mathrm{if}\; \psi|_{N_0\cap g_i Hg_i^{-1}}\; \mathrm{is\; trivial}\\ 0 &\; \mathrm{otherwise}\end{cases}.$$ 
Therefore $f$ is zero on those $g_i=ak$ such that $N_0 \cap aHa^{-1} \nsubseteq \ker(\psi)$. This  implies that the support of $f$ is contained in the set of points $g_i=ak$ with $a \in A^+(\epsilon^{-1})$ using the notation of section \ref{intro51}. 
Inside the support, the values $f(g_i)$ must satisfy 
$$ |f(g_i)| \leq (const) \delta_0^{1/2}(a)(1+\sigma_*(a))^{-r}$$ for all $r>0$.

Now we turn our attention to $\mathcal{C}^*(G)$. 
As one has the Cartan decomposition for $G$, to construct a function $u \in \schwartzg$ it suffices to set values for $u$ on $k_1a^+k_2$ where $k_1$ (resp. $k_2$) are coset representatives of the left (resp. right) $H$-cosets in $K$ and $a^+ \in A^+$. In addition to this, we require 
 $$|u(k_1 a^+k_2)|\leq (const) \Xi(a^+)(1+\sigma_*(a^+))^{-r}$$ for all $r\in \mathbb{R}$. If this bound is satisfied for all $g\in G$, then $u \in \schwartzg$. 

In order to complete the proof of Lemma \ref{pairing}(1), it suffices to prove the following statement: 

\begin{lemma}\label{schwartzsurj}
For all sufficiently small $H$, the map $(-)_{\psi}: \schwartzg \rightarrow \mathcal{C}^*(N_0\setminus G ;\psi)^{H}$ given by $$u \mapsto \int_{N_0} \psiinv(n) u(ng)\;dn$$ is surjective . 
\end{lemma}
\begin{proof}
Let $\psi$ be the nondegenerate smooth character as before. Then $$\ker(\psi) = [N_0,N_0] N_{\alpha_1, i_{\alpha_1}}\ldots N_{\alpha_r,i_{\alpha_r}}$$ 
where $\alpha_j\in \Delta.$ If we set $n= \max\{i_{\alpha_1}, \ldots, i_{\alpha_r}\}$, then $\psi$ is trivial on $N_n$. We may always assume that $H=H_m$ is sufficiently close to the identity element so that $N_m \subsetneq N_n$.

In order to prove this lemma, given any $f \in \mathcal{C}^*(N_0\setminus G;\psi)^H$, we must construct a function $u \in \mathcal{C}^*(G)^H$ which maps to $f$.

Suppose that $a \in A^+(\epsilon^{-1}) -A^+$.  Let $a^+=w^{-1} a w$ where $w \in W$ is chosen so that $a \in w A^+ w^{-1}$. Define $$u(wa^+w^{-1}k)=\frac{1}{\mathrm{meas}\;((N_0\cap aHa^{-1})N_{0,m})} f(ak).$$ Now set $u(hwa^+w^{-1} kh')=u(wa^+w^{-1}k)$ for all $h,h' \in H$. 

If $a \in A^+$, set $$u(ak)=\frac{1}{\mathrm{meas}\;N_{0,m}} f(ak).$$ As above, we require $u$ to be $H$-invariant on the right and on the left. Finally, set $u$ to be zero on every other coset. 

We claim that $u \in \schwartzg$. We have the inequality $$(const.)\delta_0^{1/2}(a) \leq \Xi(a)$$ for any $a \in A^+$. (See for instance \cite[Lemme II.1.1]{walsp}) Therefore $u$  satisfies the required growth condition on the points $a$ contained in $A^+$. 

It remains to show that if $a \in A^+(\epsilon^{-1})-A^+$, then $$\frac{\delta_0^{1/2}(a)}{\mathrm{meas}\;((N_0\cap aHa^{-1})N_{0,m})} \leq (const.)\Xi(a^+)$$ where $a \in wA^+w^{-1}$ and $a^+=w^{-1}aw.$ 

Since
$$\mathrm{meas}\;((N_0\cap aHa^{-1})N_{0,m}) =( \prod_{\substack{\alpha \in \,^{nd}\Sigma^+\\ \val(\alpha(a))<0}} \,_{\alpha}q^{-\lambda_{\alpha}\val(\alpha(a))}) \,\mathrm{meas}\;N_{0,m}$$ by definition of $N_0\cap aHa^{-1}$, 
\begin{align*}
\frac{\delta_0^{1/2}(a)}{\mathrm{meas}\;((N_0\cap aHa^{-1})N_{0,m})}
& =\prod_{\substack{\alpha \in \,^{nd}\Sigma^+ \\ \val(\alpha(a)) < 0 }} \,_{\alpha}q^{\frac{\lambda_{\alpha}}{2}\val(\alpha(a))} \prod_{\substack{\alpha \in \,^{nd}\Sigma^+ \\ \val(\alpha(a)) \geq 0 }} \,_{\alpha}q^{-\frac{\lambda_{\alpha}}{2}\val(\alpha(a))}\,(\mathrm{meas}\; N_{0,m})^{-1} \\
&=\delta_{w^{}P_{\emptyset}w^{-1}}^{1/2}(a)\, (\mathrm{meas}\; N_{0,m})^{-1}\\
&=\delta^{1/2}_0(a^+) (\mathrm{meas}\; N_{0,m})^{-1} \leq (const.) \Xi(a^+).
\end{align*}

Here $\delta_{w^{}P_{\emptyset}w^{-1}}^{1/2}$ is half the sum of all positive roots with respect to the simple root system defined by the chamber $wA^+w^{-1}$. 

Finally we show that $u_{\psi}=f$. It suffices to show that $u_{\psi}(ak)=f(ak)$ for all $a \in A^+(\epsilon^{-1})$ and $k \in K$ modulo $H_m$. By construction and using Lemma \ref{nomatchingcoset}, if $n \notin ((N_0\cap aHa^{-1}).N_{0,m})$, then $u(nak)=0$

With $a \in A^+(\epsilon^{-1})$ and $N_0\cap aHa^{-1} \nsubseteq N_{0,m}$, 

\begin{align*}u_{\psi}(a k) &= \int_{N_0} \psiinv(n) u(  na k) \; dn \\
&=\mathrm{meas} \; ((N_0 \cap aHa^{-1})N_{0,m}) \;\left[ \psiinv(1)u(a_ik) +\sum_{ \substack{n_l \in N_0/((N_0\cap aHa^{-1})N_{0,m})\\ n_l \neq 1 }} \psiinv(n_l)u(  n_l  a k)\right]\\
&=\mathrm{meas} \; ((N_0 \cap aHa^{-1})N_{0,m}). \psiinv(1) \left(\frac{f(ak)}{\mathrm{meas} \, (N_0 \cap aHa^{-1})N_{0,m}}\right) = f(ak).
\end{align*}

With $a \in A^+(\epsilon^{-1})$ and $N_0\cap aHa^{-1} \subseteq N_{0,m}$
\begin{align*}u_{\psi}(a k) &= \int_{N_0} \psiinv(n) u( n  a k) \; dn \\
&=\mathrm{meas} \, (N_0 \cap aHa^{-1}) \; \left[ \psiinv(1).\,[N_{0,m}:N_0 \cap aHa^{-1}]\,  u( ak)
+\sum_{ \substack{n_l \in N_0 / N_0 \cap aHa^{-1} \\ n_l\notin N_{0,m}}} \psiinv(n_l) u(  n_l   ak)\,\right] \\
&=\mathrm{meas}\,(N_{0,m}). \psiinv(1) \left(\frac{f(ak)}{\mathrm{meas}\,(N_{0,m})}\right) = f(ak).
\end{align*}
\end{proof}

\subsection{} 
In this section we prove Lemma \ref{pairing}(2). The terms \emph{rational}, \emph{regular}, \emph{pole} and the notion of \emph{vector valued integration} will be used in the same sense as Waldspurger defines it in \cite[Section IV-V]{walsp}.

First of all, recall the spaces  $$Wh_{\psi^{\theta}}(H_{\sigma})=\hom_{M_{\theta} \cap N_0} (\sigma, \mathbb{C}_{\psi^{\theta}})$$ and $$Wh_{\psi}(I(\sigma,\nu))= \hom_{N_0} (I_{\bar{P}_{\theta}}^{G} \sigma \boxtimes \chi_{\nu}, \mathbb{C}_{\psi}).$$ It is well known that these spaces of Whittaker functionals are one dimensional. 

For any $\lambda \in Wh_{\psi^{\theta}}(H_{\sigma})$ and  any $w \in I (\sigma)_{\tau}$ (and any $\tau$), consider the following two integral transforms:

\begin{align}
J_{\sigma, \nu}(\lambda)(w) &= \int_{{N}_{\theta}} \psi^{-1}(n) \lambda( w_{\nu}(n) )\;dn \label{jaceq1}\\
E_{\sigma, \nu}(w) &= \int_{{N}_{\theta}} \psi^{-1}(n)  w_{\nu}(n)\;dn\label{jaceq2}
\end{align} 
The following result guarantees that these transforms are well defined. 

\begin{thm}[c.f.\cite{jacq,casselman,shah}]\label{jacthes} The integral \eqref{jaceq1} (resp. \eqref{jaceq2}) converges absolutely for any $\nu \in \mathfrak{a}_{\theta, \mathbb{C}}^{*,-}$ and $w \in I(\sigma)$ and can be analytically continued to all of  $\mathfrak{a}_{\theta, \mathbb{C}}^{*}$. Then $J_{\sigma,\nu}$ defines an isomorphism of $Wh_{\psi^{\theta}}(H_{\sigma})$ onto $Wh_{\psi}(I(\sigma,\nu))$ and both $J_{\sigma,\nu}$ and $E_{\sigma,\nu}$ are smooth as functions on $\mathcal{O}$. 
\end{thm}
\begin{flushright}
$\qed$
\end{flushright}
\begin{tcor}\label{pairingconv}
Assume that $\sigma$ is a $\psi^{\theta}$-generic discrete series representation of $M_{\theta}$ and consider the induced representation $(\pi_{\bar{P}_{\theta}, \sigma,\nu}, I_{\bar{P}_{\theta}}^G\, \sigma\boxtimes \chi_{\nu} )$.  Let $\lambda$ be a Whittaker functional on $\sigma$ and  $f \in \mathcal{C}^*(N_0 \setminus G; \psi)$
The integral defined by $$ \int_{N_0 \setminus G} f(g) \overline{J_{\sigma,\nu}(\lambda)(\pi_{\bar{P}_{\theta}, \sigma, \nu}(g)w)}\, dg$$ is absolutely convergent. 
\end{tcor}
\begin{proof} This follows immediately from Theorem \ref{jacthes} and Proposition \ref{tempered}. \end{proof}
 The following lemma relates the map $J_{\sigma,\nu}$ to the Harish-Chandra transform. 
\begin{tlemma} \label{adj}
With notation as above, $$ \int_{N_0 \setminus G} \overline{J_{\sigma,\nu}(\lambda)(\pi_{\bar{P}_{\theta},\sigma,\nu}(g)w)} f(g)\; dg = \int_{(N_* \setminus M_{\theta} )\times K} \chi_{-\nu}(m)\overline{ \lambda(\sigma(m) w(k))} (R(k)f)^{P_{\theta}}(m)\;dm\,dk 
$$ for all $w \in I(\sigma)$ and all $f \in \mathcal{C}^*(N_0 \setminus G; \psi)$. 
\end{tlemma}
\begin{proof}
Recall that $N_0 = N_{\theta} N_*$. Then by definition,  \begin{align*}
\int_{N_0 \setminus G} &\overline{\lambda\left(\int_{N_{\theta}} \psiinv(n) w_{\nu}(ng)\, dn\right)} f(g)\,dg
\\
&= \int_{N_{*} \setminus G} \overline{\lambda(\w{g})} f(g) \,dg \\
&=\int_{\bar{N}_{\theta} \times (N_* \setminus M_{\theta}) \times K} \chi_{\rho-\nu}(m) \overline{\lambda(\sigma(m)w(k))} f(\bar{n}mk)\, d\bar{n}\,dm\,dk\\ \intertext{by the Iwasawa decomposition,}
&=\int_{(N_* \setminus M_{\theta}) \times K} \chi_{-\nu}(m) \overline{\lambda(\sigma(m)w(k))} \left(\int_{\bar{N}_{\theta}} \delta^{-\frac{1}{2}}_{\bar{P}_{\theta}}(m)f(\bar{n}mk)\, d\bar{n}\right)\,dm\,dk \\
&=\int_{(N_* \setminus M_{\theta}) \times K} \chi_{-\nu}(m) \overline{\lambda(\sigma(m)w(k))} (R(k)f)^{P_{\theta}}(m)\,dm\,dk .
\end{align*}
\end{proof}

\begin{proof}[Proof of Lemma \ref{pairing}(2)] Before we begin, we fix a discrete series $\sigma_0 \in \mathcal{O}$. We will use the generic symbol $\sigma$ to represent $\sigma_0\otimes \chi_{\nu}$ for any unramified character $\chi_{\nu}$ and identify the representation space $H_{\sigma_0}$ with $H_{\sigma}$. Denote the invariant Hermitian inner product on $H_{\sigma_0}$ as $\la, \ra_{\sigma_0}$. 

We may assume
$$f_{\psi}(g)= \int_{N_0} \psiinv(n_0) \int_{\mathcal{O}} \varphi_{w,v}[n_0g](\nu)\alpha(\nu)\dsig\,dn_0$$ for some datum $[\alpha, \mathcal{O}, \bar{P}_{\theta}]$ and fixed $w,v \in I(\sigma)$. For convenience, we let $g=1$ in the following discussion. 

By  Lemma \ref{gammafact}, up to a scalar $\gamma(G|M_{\theta})$, 
\begin{align*} C=f_{\psi}(1)&=\int_{N_0} \psiinv(n_0) \intorbit \left(\int_{N_1} \la w_{\nu}(n_1 n_0), v_{\nu}(n_1) \ra_{\sigma_0} \,dn_1\right)\alpha(\nu)\,\dsig \,dn_0.
\end{align*}  
where $N_1 = N_{\theta}.$ 
As the integral over $N_0\times \mathcal{O} \times N_1$ above is finite we may  switch the order of integration so that we integrate against $N_0$ followed by $N_1$ and then $\mathcal{O}$. Thus  after a change of variable, 
\newcommand{\relim}[1]{\lim_{#1\rightarrow \infty}}

$$C_{}= \intorbit \left(\int_{N_{1}\times N_{0}} \psiinv(n_0) \psi(n_1) \la w_{\nu}(n_0), v_{\nu}(n_1) \ra_{\sigma_0} \aplms dn_0\,dn_1\right)\alpha(\nu)\, \dsig. $$

We split $N_{0}= N_{1} N_{2}$ where $N_{2} = M_{\theta} \cap N_{0}$. Denote the inner integral (in brackets) over $N_1 \times N_1\times N_2$ by $I$. 

After rearranging integrals we obtain, 
\begin{align*}I&=   \int_{N_{1}}  \psi^{-1}(n_1)\int_{N_{2}} \psi^{-1} (n_2) \left< \sigma(n_2)w_{\nu}(n_1),\left(\int_{N_{1}} \psiinv(n_1) v_{\nu}(n_1)\,dn_1 \right)\right>_{\sigma_0} \aplms dn_2\,dn_1\\
&=\int_{N_1}\psiinv(n_1)\left(\int_{N_2}\psiinv(n_2)\la \sigma(n_2)w_{\nu}(n_1), E_{\sigma_0, \nu}(v)\ra_{\sigma_0}\;dn_2\right)\,dn_1.
\end{align*}\normalsize
Observe that  
$$ \la \sigma(n_2)w_{\nu}(n_1), E_{\sigma_0, \nu}(w)\ra_{\sigma_0}$$ is a matrix coefficient of $\sigma$ and hence is in $\mathcal{C}^*(M_{\theta})$ (modulo the center of $M_{\theta}$).We see from the proof of Lemma \ref{lemma42} that the integral  $$ \int_{N_2} \psiinv(n_2) \la \sigma(n_2)w_{\nu}(n_1), E_{\sigma_0, \nu}(w)\ra_{\sigma_0}$$ converges absolutely (and defines the projection to $\mathcal{C}^*(N_2\setminus M_{\theta};\psi^{\theta})$). This defines a functional $\lambda_{E_{\sigma_0, \nu}(v)} \in Wh_{\psi^{\theta}}(H_{\sigma_0})$. 

Therefore $$I_{}=\int_{N_1} \psiinv(n_1) \lambda_{E_{\sigma_0, \nu}(v)}(\w{n_1})\;dn_2.$$

By Theorem \ref{jacthes}, this integral converges giving us the linear map $J_{\sigma_0, \nu}(\lambda_{E_{\sigma_0, \nu}})(w)$. Thus we have
$$C= \intorbit J_{\sigma_0, \nu}(\lambda_{E_{\sigma_0, \nu}(v)})(w) \alpha(\nu)\, \dsig.$$ 
 
With this we can complete the proof. Now we obtain $$f_{\psi}(g)=\intorbit J_{\sigma_0, \nu}(\lambda_{E_{\sigma_0, \nu}(v)})(\pi_{\bar{P}_{\theta}, \sigma}(g)w)\alpha(\nu)\;\dsig$$ by replacing $w$ with $\pi_{\bar{P}_{\theta}, \sigma}(g)w$ in the discussion above. 

Then
\begin{align*}B(f,\varphi)&=\int_{N_0\setminus G} \overline{\varphi(g)} \intorbit J_{\sigma_0, \nu}(\lambda_{E_{\sigma_0, \nu}(v)})(\pi_{\bar{P}_{\theta}, \sigma}(g)w)\alpha(\nu)\;\dsig\,dg.
\end{align*} The interchange of integrals is permitted since the integral defining $B$ converges absolutely by Lemma \ref{lemma42}. Now we apply Lemma \ref{adj} to see that $B(f, \varphi)$ is 
\begin{align}
\intorbit \left(\int_{N_2 \setminus M_{\theta} \times K} \chi_{\nu}(m)\lambda_{E_{\sigma_0, \nu}(v)}(\sigma_0(m)w(k))\overline{(R(k)\varphi)^{P_{\theta}}(m)}\, dm\,dk\right)\, \alpha(v) \,\dsig. \label{reff5} 
\end{align}
Thus if $f \in\mathcal{C}^*_{cont}(G)$, $B(f, \varphi)=0$ for all $\varphi \in \,^{\circ}\mathcal{C}^*(N_0 \setminus G;\psi)$ . 

\end{proof}


\section{The Plancherel formula for $L^2(N_0 \setminus G ; \psi)$}

\subsection{}\label{sec61}
Suppose $P_{\theta}=M_{\theta}N_{\theta}$ is a standard parabolic subgroup of a quasi-split group $G$ and $(\sigma, H_{\sigma})$ a discrete series representation of $M_{\theta}$. Assume that $\sigma$ is $\psi^{\theta}$ generic and let $Wh_{\psi^{\theta}}(H_{\sigma})$ be the space of its Whittaker functionals. As $M_{\theta}$ is also quasi-split, a result of Shalika assures us that this space is one-dimensional. 

We recall the following construction used in the proof of Lemma \ref{pairing}(2). As before let $N_* = M_{\theta} \cap N_0$. Then $N_0 = N_{\theta} N_*$. Let $w, v \in H_{\sigma}$ and consider the transform: 
$$\lambda_v(w):=\int_{N_*} \psiinv(n_*)\la \sigma(n_*)w, v \ra_{\sigma}\, dn_*. $$

This is well defined as $\la \sigma(m)w,v\ra_{\sigma} \in \mathcal{C}^*(M_{\theta}).$ The integral above converges absolutely for any $w \in H_{\sigma}$ and defines a Whittaker functional $\lambda_v \in Wh_{\psi^{\theta}}(H_{\sigma})$. Observe that $\lambda_v(w)= \overline{\lambda_w(v)}$. 

Choose any $ \lambda \in Wh_{\psi^{\theta}}(H_{\sigma})$ and $w \in H_{\sigma}$ satisfying $\lambda(w) = 1.$ Then for $u \in H_{\sigma}$, we have $$\lambda_u = \lambda_u(w)\lambda =  \overline{\lambda_{w}(u)} \lambda.$$

\begin{lemma} \label{whitfunctbasis}
Let $\alpha \in C^{\infty}(\mathcal{O})$, $w,v \in I(\sigma).$ We define $\eta=\lambda_w$.  Then $$(f_{[\alpha, \mathcal{O}, \bar{P}_{\theta}];w,v})_{\psi}(g)= \int_{\mathcal{O}} \overline{J_{\sigma, \nu}(\eta)(w)}J_{\sigma, \nu}(\lambda)(\pi_{\bar{P}_{\sigma}, \sigma, \nu}(g)v) \alpha(\nu) \,\dsig.$$
\end{lemma}
\begin{proof}
We already know that $$(f_{[\alpha, \mathcal{O}, \bar{P}_{\theta}]; w,v})_{\psi}(g)= \int_{\mathcal{O}} J_{\sigma, \nu}(\lambda_{E_{\sigma, \nu}(w)})(\pi_{\bar{P}_{\theta}, \sigma, \nu}(g)v)\alpha(\nu)\,\dsig$$ from the proof of Lemma \ref{pairing}(2). 

However we may write $\lambda_{E_{\sigma, \nu}(w)}= \overline{\eta(E_{\sigma, \nu}(w))}\lambda =  \overline{J_{\sigma, \nu}(\eta)(w)} \lambda$ from which the lemma follows immediately. 
\end{proof}

Now fix a base point $\sigma_0 \in \mathcal{O}$ and for every $\sigma \in \mathcal{O}$ identify $H_{\sigma}=H_{\sigma_0}$. Then without loss of generality, $Wh_{\psi^{\theta}}(H_{\sigma_0})=Wh_{\psi^{\theta}}(H_{\sigma})$. Thus, there is an unambiguous choice of basis $\lambda \in Wh_{\psi^{\theta}}(H_{\sigma})$ for each $\sigma \in \mathcal{O}$. The space $C^{\infty}(\mathcal{O}; Wh_{\psi^{\theta}}(H_{\sigma}))$ will then consist of all functions of the form $$\{(\nu \mapsto \alpha(\nu)\lambda)\mid \alpha \in \mathcal{C}^{\infty}(\mathcal{O})\}.$$

Define a map $\Lambda: C^{\infty}(\mathcal{O}) \otimes I_{\sigma}\rightarrow C^{\infty}(\mathcal{O}; Wh_{\psi^{\theta}}(H_{\sigma}))$ by the formula $$\Lambda(\alpha \otimes f)(\nu) = \alpha(\nu) \overline{J_{\sigma, \nu}(\eta)(f)}\lambda.$$
\begin{lemma}\label{surjection1}
The linear map $\Lambda$ defined above is surjective.
\end{lemma}
\begin{proof}
We observe that $\lambda_w(v)=(\varphi_{v,w})_{\psi^{\theta}}(1).$ Thus we may as well assume $\eta$ is not the zero functional. 

We must prove this assertion for $(\nu \mapsto \alpha(\nu)\lambda)$ for any $\alpha \in C^{\infty}(\mathcal{O})$ and a fixed $\lambda \in Wh_{\psi^{\theta}}(H_{\sigma})$. 

By Theorem \ref{jacthes} $J_{\sigma, \nu}(\eta_j)$ is a nonzero functional. Thus one may find $w \in I(\sigma)$ such that $J_{\sigma, \nu}(\eta)(w_{\nu})\neq 0.$ 
As $J_{\sigma, \nu}$ is continuous in $\nu$, $J_{\sigma, \mu}(\eta)(w_{\mu})$ is nonzero as long as we vary $\mu$ around a small enough neighborhood $U_{\nu}$ of $\nu \in \mathcal{O}$. 

Since $\mathcal{O}$ is compact, we may find a finite covering $\{U_l\}_{l=1}^p$ of the support of $\alpha \in C^{\infty}(\mathcal{O})$ subject to the condition that for each of these open sets $U_l$, there exists vectors $w^l \in I(\sigma)$ such that on each of these open neighborhoods $U_l$, $$J^l_{\sigma,\mu}=J_{\sigma, \mu^{}}(\eta)(w^l_{\mu^{}})_{}$$ is nonzero for all $\mu^{}\in U_l.$

Now choose a partition of unity, $\{\varphi_l\}$ subordinate to the cover $\{U_l\}$  and define for each open set $U_l$ a  (smooth) function $g_l(\nu)$ which vanishes outside $U_l$ and is equals to $(g_{}^l(\nu))_{}= \varphi_l(\nu)\overline{(J^{l}_{\sigma, \nu})^{-1}}$.

Thus we clearly have \begin{align*}
\sum_{l}\Lambda(g_l \alpha \otimes w^l_{\nu}) &= \sum_{l}g_l(\nu)\alpha(\nu)\overline{J_{\sigma, \nu}(\eta)(w^l_{\nu})}\lambda \\ 
&=\left(\sum_l\varphi_l(\nu) \right)\alpha(\nu) \lambda_{}=\alpha(\nu)\lambda.\end{align*} This proves the lemma. 
\end{proof}
Consider the following two Hermitian inner products on the discrete series representation $H_{\sigma}$ and the image of $\sigma$ in $\,^{\circ} \mathcal{C}^*(N_2 \setminus M_{\theta}; \psi^{\theta}):$ $ \la v, w \ra_{\sigma}$ and \begin{align}\int_{N_2 \setminus M_{\theta}} \lambda(\sigma(m)v) \overline{ \eta(\sigma(m) w)}\, dm\label{reff4}\end{align} for $\lambda , \eta \in Wh_{\psi^{\theta}}(H_{\sigma})$ respectively. By Schur's lemma, the invariant Hermitian form on an irreducible unitary representation is unique up to scalar. Thus, if we set $(\lambda, \eta)_{\sigma}$ to be the ratio of \eqref{reff4} and $\la v,w \ra_{\sigma}$, this pairing is sesquilinear and positive definite.

\begin{lemma}\label{bilinear identity} 
For any two vectors $w, v \in H_{\sigma}$, $$d(\sigma)(\lambda_v, \lambda_w)_{\sigma}=\lambda_v(w)$$ where $d(\sigma)$ is the formal degree of $\sigma$.
\end{lemma} 
\begin{proof}
Let $M_{\theta}=M$ and by abuse of notation write $\psi^{\theta}$ as $\psi$. 
By definition \begin{align*} (\lambda_v, \lambda_w)_{\sigma} \la x,y \ra_{\sigma}&=\int_{N_{*} \setminus M} \lambda_v(\sigma(m)x) \overline{ \lambda_w(\sigma(m)y)}\, dn \\
&= \int_{N_{*} \setminus M} \int_{N_{*} \times N_{*}}  \psiinv(n_1n_2^{-1})\la \sigma(n_1m)x, v \ra_{\sigma} \overline{ \la \sigma(n_2m)y,w \ra}_{\sigma}\, dn_1dn_2 dm. \\
\intertext{By change of variable $n_1n_2$ to $n_1$, this equals}
&\;\;\;\,\int_{N_{*} \setminus M} \int_{N_{*} \times N_{*}} \psiinv(n_1) \la \sigma(n_2m)x, \sigma(n_1)^{-1} v \ra_{\sigma} \overline{ \la \sigma(n_2m)y,w\ra}_{\sigma}\,dn_1 dn_2 dm \\
&= \int_M \int_{N_{*}} \psiinv(n) \la \sigma(m)x, \sigma(n)^{-1}v\ra_{\sigma} \overline{\la\sigma(m)y,w\ra}_{\sigma} \, dn  dm. \\
\intertext{By switching the order of integration to integrate over $M$ first, we get}
&\;\; \;\,d(\sigma)^{-1} \la x,y \ra_{\sigma}  \int_{N_{*}} \psiinv(n)\la \sigma(n)w, v \ra_{\sigma}\,dn.\\
&=d(\sigma)^{-1} \la x,y\ra_{\sigma} \lambda_v(w).
\end{align*}
By cancelling away $\la x,y \ra_{\sigma}$, we obtain our stated identity. 
\end{proof}

\begin{lemma} \label{lastindentity}
Let $\beta' \in C^{\infty}(\mathcal{O}; Wh_{\psi^{\theta}}(H_{\sigma}))$, $\alpha \in C^{\infty}(\mathcal{O})$ and $v \in I(\sigma)$. Then 
$$(\Lambda(\alpha\otimes v),\beta')_{\sigma} =d(\sigma)^{-1}\alpha(\nu) \overline{\beta(\nu)J_{\sigma,\nu}(\lambda)(v)}$$ where $\beta \in C^{\infty}(\mathcal{O})$ defines $\beta'$ (i.e. $\beta'(\nu)=\beta(\nu)\lambda$). 
\end{lemma}
\begin{proof}
To justify this identity, we need only to observe that by definition, $$(\Lambda(\alpha\otimes v),\beta')_{\sigma} = \alpha(\nu) \overline{\beta(\nu)J_{\sigma,\nu'}(\eta')(v)}(\lambda', \lambda)_{\sigma}=\alpha(\nu) \overline{\beta(\nu)J_{\sigma,\nu'}(\eta'(\lambda, \lambda')_{\sigma})(v)}$$ where $\lambda'$ is a choice of basis for $Wh_{\psi^{\theta}}(H_{\sigma})$ satisfying $(\lambda',\lambda')_{\sigma}=1$ and $\eta'=\lambda_{w'}$ where $w' \in H_{\sigma}$ is chosen so that $\lambda'(w')=1$.

 However, we note that $\lambda_v = \overline{\eta'(v)}\lambda' = d(\sigma) (\lambda_v, \eta')\lambda'$ for any $v \in I(\sigma)$ by Lemma \ref{bilinear identity}. In other words $\eta'= d(\sigma)^{-1} \lambda'$. This implies the result. 
\end{proof}

\subsection{}\label{sec62}
We give a quick overview of the theory of intertwining operators for induced representations adequate for our purposes. Recall that $W$ denotes the Weyl group of $G$. More precisely, we write this as $W^G$ for if $M$ is a Levi subgroup of a standard parabolic subgroup $P_{\theta}$, $W^M$ will denote the Weyl group of $M$. 

Now assume that we are given two 
 standard parabolic subgroups, $P$ and $P'$ with Langlands decomposition $P=MN$ and $P=MN'$ respectively. For any $w \in W(G|M)$, let $k_w$ be a choice of representative of $w$ in $G$. Then define $w.P=k_w P k_w^{-1}$ and if $\sigma$ is a representation of the Levi $M$,  define $$w\sigma(m):=\sigma(k_w^{-1}mk_w)$$  for $m \in w.M$. In addition to this, define $$W(G|M):=\{w \in W^G\mid w.M = M \}/W^M$$
 
If $P, P'$ are two standard parabolic subgroups with Langlands decomposition $P=MN$ and $P=MN'$ respectively, then define  
$$J_{P|P'}(\sigma)f(g) = \int_{N \cap N'\setminus N'} f(n'g)\, dn'$$ with $f \in I_P^G \sigma$ provided this integral converges. It is well known that $J_{P|P'}(\sigma)$ can be meromorphically continued to a rational function on $\mathcal{O}_{\mathbb{C}}$. (See \cite[Theoreme IV.1.1. and pg. 276-278]{walsp} for nomenclature and the precise statement.) 

Now we assume $\sigma$ is unitary. Recall the unitary intertwining operator defined in \cite[pg. 295]{walsp} ($ w \in W(G|M)$ and $ \sigma \in \mathcal{O}$) $$\,^{\circ}c_{P'\mid P}(w, \sigma): I_{P}^{G} \sigma \otimes I_{P}^G \sigma^{\vee}\rightarrow I_{P'}^G\, w\sigma \otimes I_{P'}^G w\sigma^{\vee}$$ for two parabolic subgroups $P$ and $P'$ with the same Levi subgroup $M$. 

It is known that $\,^{\circ}c_{P'\mid P}(w, \sigma)$ is regular on $\mathcal{O}$ and that we may express this map in the following form: $$u\otimes v \mapsto J_{P'|w.P}(w\sigma)L(k_w)u \otimes J_{w.P|P'}(w\sigma^{\vee})^{-1}L(k_w)v$$ where $L$ denotes left translation. See \cite[proof of Lemme V.3.1]{walsp}. We define $$A_w(\sigma):=J_{P'|w.P}(w\sigma)L(w)$$
and for a fixed $\sigma \in \mathcal{O}$, consider $A_w(\sigma\boxtimes \chi_v)$ as a function on $\mathcal{O}$ in which case we write $A_w(\sigma \boxtimes \chi_{\nu})$ as $A_w(\nu)$.

Given $\alpha \in C^{\infty}(\mathcal{O}; Wh_{\psi^{\theta}}(H_{\sigma}))$ and $v \in I(\sigma)$ we define
\begin{align}W_{[\alpha, \mathcal{O}, \bar{P}_{\theta}]; v}(g):= \int_{\mathcal{O}} J_{\sigma, \nu}(\alpha(\nu))(\pi_{\bar{P}_{\theta}, \sigma, \nu}(g)v) \,\dsig.\label{whittaker}\end{align}
\newcommand{\whittaker}{W_{[\alpha, \mathcal{O}, \bar{P}_{\theta}]; v}}
Let $\beta' \in  C^{\infty}(\mathcal{O}; Wh_{\psi^{\theta}}(H_{\sigma}))$ be defined by $\beta \in C^{\infty}(\mathcal{O})$.

\begin{lemma}\label{lemma621} Let $A_{w}(\nu)$ be the intertwining operator defined by the unitary map $\,^{\circ}c_{\bar{P_{\theta}}|\bar{P_{\theta}}}(w,\sigma)$. Then
there exists a map $$M_w(\nu)^{-1}:Wh_{\psi^{\theta}}(H_{\sigma}) \rightarrow Wh_{\psi^{\theta}}(H_{w\sigma})$$ smooth on $\mathcal{O}$ so that, 
\begin{align}W_{[\beta', \mathcal{O}, \bar{P}_{\theta}]; v}(g)= \int_{\mathcal{O}} \beta(\nu)J_{w\sigma, w\nu}(M_{w}(\nu)^{-1}\lambda)(\pi_{\bar{P}_{\theta}, w\sigma, w\nu}(g)A_w(\nu)v)\, \dsig\label{reff14}\end{align} for any $w \in W(G|M_{\theta})$. 
\end{lemma}
\begin{proof}
Indeed, if $\lambda \in Wh_{\psi^{\theta}}(H_{w\sigma})$, then $J_{w\sigma, w\nu}(\lambda) \circ A_w(\nu) \in J_{\sigma, \nu}(Wh_{\psi^{\theta}}(H_{\sigma}))$. 

Let $Z$ denote the poles of $A_w(\nu)$ on $\mathcal{O}$ and let $\mathcal{O}'=\mathcal{O}-Z$. Then
as $J_{\sigma, \nu}$ is an isomorphism (Theorem \ref{jacthes}), we define $M_w(\nu)$ by 
$$J_{w\sigma, w\nu}(\lambda)\circ A_w(\nu) = J_{\sigma, \nu}(M_w(\nu)\lambda)$$ for all $\nu \in \mathcal{O}'$. Since $A_w(\nu)$ is an isomorphism on the points where it is regular, $M_w(\nu)^{-1}$ is well defined on $\mathcal{O}'$. 

If $\lambda \in Wh_{\psi^{\theta}}(H_{\sigma})$, then \begin{align} J_{w\sigma, w\nu}(M_w(\nu)^{-1}\lambda) \circ A_w(\nu) = J_{\sigma, \nu}(\lambda)\end{align} on $\mathcal{O}'$. Since the right hand side is defined on $\mathcal{O}$, we may extend this equality to all of $\mathcal{O}$ by setting $M_w(\nu)^{-1}=0$. The claim follows immediately from definition. 
\end{proof}

We remark that if $P=G$, then we take $\mathcal{O}=\{e\}$ as a singleton so that $\alpha$ is identified with $\lambda \in Wh_{\psi}(H_{\sigma})$ where $\sigma$ is a discrete series representation of $G$. Then $W_{[\alpha,\{ e\},G];v}(g)$ in \eqref{whittaker} is $W(v,\lambda)(g)$. 

\begin{lthm}\label{main4}

\begin{enumerate} 
\item If $\alpha \in C^{\infty}(\mathcal{O}; Wh_{\psi^{\theta}}(H_{\sigma}))$, then $\whittaker \in \mathcal{C}^*(N_0 \setminus G ; \psi)$ for all $v \in I(\sigma)$. 
\item 
The span of $\whittaker$  over all $v \in I(\sigma) $ and datum $[\alpha, \mathcal{O}, \bar{P}_{\theta}]$  is equal to  $C^*(N_0 \setminus G; \psi)$. In particular, the span is dense in $L^2(N_0 \setminus G; \psi)$. 
\item 
Given two datum $[\alpha, \mathcal{O}, \bar{P}_{\theta}]$ and $[\beta, \mathcal{P}, \bar{P}_{\theta'}]$, suppose either $\theta\neq \theta'$ or  $\mathcal{O} \neq w\mathcal{P}$ for any $w \in W(G|M_{\theta})$, then $\la \whittaker, W_{[\beta, \mathcal{P}, \bar{P}_{\theta'}];w}\ra=0.$ 

Otherwise, if $\theta = \theta'$ and $\mathcal{O}= \mathcal{P}$, then $$\la \whittaker,  W_{[\beta, \mathcal{O}, \bar{P}_{\theta}];w} \ra = |W(G|M_{\theta})|^{}\gamma(G|M_{\theta})^{} c(G|M_{\theta})^2 \left(\int_{\mathcal{O}} (\alpha(\nu), \beta(\nu))_{\sigma} \dsig\,\right) \la v,w\ra_{I(\sigma)}.$$ 
\end{enumerate} 
\end{lthm}

\begin{proof}[Proof of Theorem \ref{main4} (1)]
This is an immediate consequence of Lemma \ref{whitfunctbasis} and Lemma \ref{surjection1}. One first expresses $\alpha \in C^{\infty}(\mathcal{O}; Wh_{\psi^{\theta}}(H_{\sigma}))$ as a finite combination of $\Lambda(\beta \otimes w).$ By construction, it is immediate that $\whittaker$ is a finite combination of $(f_{[\beta, \mathcal{O}, \bar{P}_{\theta}];w,v})_{\psi}$. This proves the assertion.  
\end{proof}
\begin{proof}[Proof of Theorem \ref{main4} (2)] We know from the main theorem of the Plancherel formula for $L^2(G)$ that $\mathcal{C}^*(G)$ is the span over all  $f=f_{[\alpha, \mathcal{O}, \bar{P}_{\theta}]; w,v}$ for all datum $[\alpha, \mathcal{O}, \bar{P}_{\theta}]$. (When $P=G$, $f$ reduces to matrix coefficients of a discrete series representation of $G$).  We also know that taking $\psi$-coinvariants of $\mathcal{C}^*(G)$ is surjective onto $\mathcal{C}^*(N_0 \setminus G;\psi)$ by Lemma \ref{pairing} (1). 

Now let $$ \beta(\nu)=\alpha(\nu)\overline{J_{\sigma,\nu}(\eta)(w)}$$ where $\alpha \in C^{\infty}(\mathcal{O})$, $w \in I(\sigma)$ and $\eta$ defined as in section \ref{sec61}. We know that $J_{\sigma, \nu}$ is smooth as a function of $\nu$ (c.f. Theorem \ref{jacthes}), thus $\overline{J_{\sigma,\nu}(\eta)(w)}$ is smooth too. Thus $\beta \in C^{\infty}(\mathcal{O})$. Let $\beta'(\nu)=\beta(\nu)\lambda$. Then by Lemma \ref{whitfunctbasis}, it is clear that $$(f_{[\alpha, \mathcal{O}, \bar{P}_{\theta}];w,v})_{\psi} = W_{[\beta', \mathcal{O}, \bar{P}_{\theta}];v}.$$
This proves (2). 
\end{proof}


\subsection{}\label{sec63}

In this section we present the proof of Theorem \ref{main4} (3).

Given two standard parabolic subgroups of $G$, $P_{\theta}$ and $P_{\vartheta}$ we say that $(P_{\theta}, A_{\theta})$ \emph{dominates} $(P_{\vartheta}, A_{\vartheta})$ if and only if $P_{\theta} \supset P_{\vartheta}$  and $A_{\theta} \subset A_{\vartheta}$.  

\begin{lemma}\label{parabolic dominance}
Consider two standard parabolic pairs $(P_{\vartheta}, A_{\vartheta})$ and $(P_{\theta}, A_{\theta})$, a datum $[ \alpha, \mathcal{O}, \bar{P}_{\theta}]$ and suppose that $((f_{[\alpha, \mathcal{O}, \bar{P}_{\theta}];w,v})_{\psi})^{P_{\vartheta}}\neq 0. $ Then $(P_{\vartheta}, A_{\vartheta})$ dominates $(P_{\theta}, A_{\theta})$. 
\end{lemma}
\begin{proof}
Recall that $$f(g)=f_{[\alpha, \mathcal{O}, \bar{P}_{\theta}];v,w}(g)=\int_{\mathcal{O}}\la \pi_{\bar{P}_{\theta}, \sigma, \nu}(g)v_{\nu}, w_{\nu}\ra \alpha(\nu)\, \dsig.$$ By Lemma \ref{switching order} we may first integrate over $\bar{N}_{\vartheta}$ in the expression of $(f_{\psi})^{P_{\theta}}$ as the following  integral over $\bar{N}_{\vartheta} \times N_0$ 
$$ \delta_{P_{\vartheta}}^{\frac{1}{2}}(m) \int_{\bar{N}_{\vartheta}} \int_{N_0} \psiinv(n_0) f(n_0 \bar{n}m) dn_0d\bar{n}.$$

Now decompose $n_0 = k(n_0) m(n_0) \bar{n}(n_0)$ according to the Iwasawa decomposition, $G=K\bar{P}_{\vartheta}=KM_{\vartheta} \bar{N}_{\vartheta}$. Then\begin{align} \int_{\bar{N}_{\vartheta}} f(n_0 \bar{n}m)\, d\bar{n} &= \int_{\bar{N}_{\vartheta}} \int_{\mathcal{O}} \la \pi_{\bar{P}_{\theta}, \sigma, \nu}(k(n_0)m(n_0)\bar{n}m)v_{\nu}, w_{\nu}\ra \alpha(\nu)\, d\mu(\nu)\notag\\
&= \delta_{P_{\vartheta}}^{}(m(n_0))\int_{\bar{N}_{\vartheta}}\int_{\mathcal{O}}  \varphi_{v, \pi(k(n_0))w}[\bar{n}m(n_0) m](\nu) \alpha(\nu)\, d\mu(\nu)\notag\\
&= \delta_{P_{\vartheta}}^{\frac{1}{2}}(m(n_0)) \delta^{-\frac{1}{2}}_{P_{\vartheta}}(m)(f_{[\alpha, \mathcal{O}, \bar{P}_{\theta}]; v, \pi(k(n_0))w})^{P_{\vartheta}}(m(n_0)m). \label{reff10}
\end{align}
 To ensure that the right hand side is not equivalently zero, it is necessary that $(P_{\vartheta}, A_{\vartheta})$ dominates $(P_{\theta}, A_{\theta})$ by \cite[Proposition VI.4.1]{walsp} .
\end{proof}

We recall that $\la, \ra $ is the $L^2$ inner product on the space $L^2(N_0 \setminus G ; \psi)$. 

\begin{lemma}\label{lemmaorth} If two parabolic subgroups, $P_{\theta}$ and $P_{\vartheta}$ are not equal, then for any datum $[\beta, \mathcal{P}, \bar{P}_{\vartheta}]$, $[\alpha, \mathcal{O}, \bar{P}_{\theta}]$ and any discrete series representation $\varrho \in \mathcal{P}$ and $\sigma \in \mathcal{O}$, 
 $$ \la ( f_{[\beta, \mathcal{O}_Q, \bar{Q}]; w,v})_{\psi}, (f_{[\alpha, \mathcal{O}, \bar{P}]; w',v'})_{\psi} \ra =0$$ for all $w,v \in I(\varrho)$ and $ w',v' \in I(\sigma)$. 
\end{lemma}
\begin{proof}
Set $f=  f_{[\beta, \mathcal{P}, \bar{P}_{\vartheta}]; w,v}$ and $\varphi = (f_{[\alpha, \mathcal{O}, \bar{P}_{\theta}]; w',v'})_{\psi} $. Then $B(f, \varphi)$ is equals to the inner product we wish to compute above. However we know from \eqref{reff5} and  Lemma \ref{parabolic dominance} that $B(f, \varphi)=0$ unless $(P_{\vartheta},A_{\vartheta})$ dominates $(P_{\theta}, A_{\theta})$. Reversing the roles of $f$ and $\varphi$, we also conclude that $B(f, \varphi)=0$ unless $(P_{\theta},A_{\theta})$ dominates $(P_{\vartheta}, A_{\vartheta})$. Combining these two statements gives the lemma. 
\end{proof}

We wish to compute $((f_{[\alpha, \mathcal{O},\bar{P}];w,v})_{\psi})^P.$ Recall from \cite[Proposition VI.4.1]{walsp}  that $(f_{[\alpha, \mathcal{O},\bar{P}];w,v})^P$ is equals to \begin{align}
\gamma(G|M) c(G|M)^2 \int_{\mathcal{O}} \alpha(\nu)\sum_{ w \in W(G|M)}  \la \chi_{w\nu}\otimes w\sigma(m).A_w(\nu)u_{\nu}(1), A_w(\nu)v_{\nu}(1)\ra_{w\sigma}\, d\nu \label{reff11}\end{align} where $A_w(\nu)$ are the intertwining operators defining the unitary intertwining map $\,^{\circ}c_{\bar{P}\mid \bar{P}}(w, \chi_{\nu} \otimes \sigma)$. We remind the reader that $\la , \ra_{w\sigma}$ denotes the $G$-invariant \emph{Hermitian} inner product on $w\sigma$.

\begin{lprop} The Harish-Chandra transform, \begin{align} ((&f_{[\alpha, \mathcal{O},\bar{P}];u,v})_{\psi})^P \notag\\
&=\gamma(G | M) c(G|M)^2 \int_{\mathcal{O}} \alpha(\nu)\sum_{w \in W(G|M)} \chi_{w \nu}(m) 
 \lambda_{E_{w\sigma, w\nu}(A_w(\nu)v)}(w\sigma(m)A_w(\nu)u_{\nu}(1))\, d\nu.\label{reff12} \end{align}
 \end{lprop}
 \begin{proof} 
Combining the expressions in the right hand side of  \eqref{reff10}  with \eqref{reff11}, we have that \begin{align*}
&((f_{[\alpha, \mathcal{O}, \bar{P}];u,v})_{\psi})^{P}\\ &=\delta^{\frac{1}{2}}_P(m)\int_{\bar{N}}\int_{N_0}\psiinv(n_0)f(n_0\bar{n}m)\, dn_0d\bar{n}. \\
\intertext{By \eqref{reff11}, this equals to}
&=\gamma(G|M) c(G|M)^2\int_{N_0}\psiinv(n_0) \\
&\quad\quad\quad\delta^{\frac{1}{2}}_{P}(m(n_0))\int_{\mathcal{O}} \alpha(\nu)\sum_{ w \in W(G|M)}\chi_{w\nu}(m(n_0)m) \la w\sigma(m(n_0)m).A_w(\nu)u_{\nu}(1), A_w(\nu)v_{\nu}(k(n_0)^{-1})\ra_{w\sigma}\, d\nu. \\
\intertext{By moving $\delta^{\frac{1}{2}}_P(m(n_0)), w\sigma(m(n_0))$ and $\chi_{w\nu}(m(n_0))$ into the second variable of $\la, \ra_{w\sigma}$, we obtain}
& \gamma(G|M) c(G| M)^2\int_{N_0}\psiinv(n_0) \\
&\quad\quad\quad\int_{\mathcal{O}} \alpha(\nu)\sum_{ w \in W(G|M)} \chi_{w\nu}(m)  \la w\sigma(m).A_w(\nu)u_{\nu}(1), \delta_{\bar{P}}^{\frac{1}{2}}\chi_{w\nu}\otimes w\sigma(m(n_0)^{-1}).A_w(\nu)v_{\nu}(k(n_0)^{-1})\ra_{w\sigma}\, d\nu \\
&=\gamma(G| M) c(G|M)^2\int_{N_0}\psiinv(n_0) \\
&\quad\quad\quad\int_{\mathcal{O}} \alpha(\nu)\sum_{ w \in W(G|M)} \chi_{w\nu}(m)  \la w\sigma(m).A_w(\nu)u_{\nu}(1), A_w(\nu)v_{\nu}((\bar{n}(n_0)m(n_0)k(n_0))^{-1})\ra_{w\sigma}\, d\nu \\
\intertext{since $A_w(\nu) v_{\nu} \in I_{\bar{P}}^G \;w\sigma \boxtimes \chi_{w\nu}$. Continuing, this expression}
&=\gamma(G| M) c(G| M)^2 \int_{\mathcal{O}} \alpha(\nu)\sum_{ w \in W(G|M)}\chi_{w\nu}(m)  \\
&\quad\quad\quad\int_{N_0}\psiinv(n_0) \la  w\sigma(m).A_w(\nu)u_{\nu}(1),A_w(\nu)v_{\nu}(n_0^{-1}) \ra_{w\sigma}\, d\nu 
\end{align*} 
Arguing as in the proof of Lemma \ref{pairing} (2), this expression is equal to 
\begin{align*}\gamma(G | M) c(G|M)^2 \int_{\mathcal{O}} \alpha(\nu)\sum_{w \in W(G|M)} \chi_{w \nu}(m) 
 \lambda_{E_{w\sigma, w\nu}(A_w(\nu)v)}(w\sigma(m)A_w(\nu)u_{\nu}(1))\, d\nu. 
\end{align*}
\end{proof}
\begin{proof} [Proof of Theorem \ref{main4} (3)] 
Assume $P_{\theta} \neq P_{\vartheta}$. As we may express $W_{[\alpha, \mathcal{O}, \bar{P}]; v}$ in terms of $f_{\psi}$ (c.f. proof of Theorem \ref{main4} (1)), the orthogonality statement follows immediately from Lemma \ref{lemmaorth}.  

Now assume that $\theta = \vartheta$ and that we have two orbits $\mathcal{O}$ and $\mathcal{P}$ containing discrete series representations $\sigma$ and $\varrho$ respectively. Consider functions $W_{[\alpha, \mathcal{O}, \bar{P}_{\theta}]; v}$  and $W_{[\beta, \mathcal{P}, \bar{P}_{\theta}];v_1}$ where $\alpha=\Lambda(\alpha' \otimes u)$ and $\beta=\Lambda(\beta' \otimes u_1)$ with $u,v \in I(\sigma)$ and $u_1,v_1 \in I(\varrho)$. We observe that by Lemma \ref{whitfunctbasis}, $$W_{[\Lambda(\alpha'\otimes v), \mathcal{O}, \bar{P}_{\theta}];u}(g) = (f_{[\alpha, \mathcal{O}, \bar{P}_{\theta}]; u, v})_{\psi}(g).$$ 

Since Lemma \ref{surjection1} tells us that $\Lambda$ is surjective what must be shown is that if  $ \mathcal{P}\nsubseteq W(G|M_{\theta}) \mathcal{O}$, then $$ \la W_{[\beta', \mathcal{P}, \bar{P}_{\theta}];v_1},W_{[\alpha', \mathcal{O}, \bar{P}_{\theta}];v}  \ra =0.$$ To this end, it suffices to show that $$B( f_{[\beta, \mathcal{P}, \bar{P}_{\theta}]; u_1, v_1},(f_{[\alpha, \mathcal{O}, \bar{P}_{\theta}];u,v})_{\psi} )=0. $$ 

We begin from the expression of \eqref{reff5} and write the expression $\lambda_{E_{\varrho, \nu}(v_1)}(\varrho(m)u_1(k))$ there as $$W(u_1(k), \lambda_{E_{\varrho, \nu}(v_1)})(m)=:W^{u_1(k),\lambda_{Ev_1}}_{\varrho, \nu}(m).$$ In the same way $\lambda_{E_{w\sigma, w\nu}(A_w(\nu)v)}(w\sigma(m)A_w(\nu)u_{\nu}(1))$ in \eqref{reff12} is written as $$W(A_w(\nu)u_{\nu}(1),\lambda_{E_{w\sigma, w\nu}(A_w(\nu)v)})(m)=:W^{Au(1),\lambda_{EAv}}_{w\sigma, w\nu}(m).$$ If we fix a $\chi_{\nu} \in \mathrm{Im}\, X_{ur}(M_{\theta})$, then these are square integrable Whittaker functions on the Levi, $M_{\theta}$ and realizes the representations $\varrho$ and $w\sigma$ in $\,^{\circ}\mathcal{C}(N_{*}\setminus M_{\theta}; \psi^{\theta})$. 

Now  we let $\varphi=(f_{[\alpha, \mathcal{O}, \bar{P}_{\theta}];u,v})_{\psi}$ and $f= f_{[\beta, \mathcal{P}, \bar{P}_{\theta}]; u_1, v_1}$. If we substitute this into \eqref{reff5}, and use \eqref{reff12} then the expression becomes
\begin{align} 
\gamma (G| M_{\theta})&c(G| M_{\theta})^2 \int_{N_{*} \setminus M_{\theta} \times K}\int_{\mathcal{P}}\beta(\nu_1) \chi_{\nu_1}(m)W^{u_1(k),\lambda_{Ev_1}}_{\varrho, \nu_1}(m)\,d\mu_{\mathcal{P}}(\nu_1) \notag\\
&\overline{\left(\
\int_{\mathcal{O}} \alpha(\nu) \sum_{w \in W(G|M_{\theta})} \chi_{w\nu}(m)W^{Au(k),\lambda_{E Av}}_{w\sigma, w\nu}(m)\, d\nu \right)}\,dm\,dk \notag \\
=\gamma (G| M_{\theta})&c(G| M_{\theta})^2\int_{\mathcal{P}}\beta(\nu_1)\int_{\mathcal{O}}\overline{\alpha(\nu)} \\
&\int_{N_{*} \setminus M_{\theta} \times K} \sum_{w \in W(G|M_{\theta})}\chi_{\nu_1}(m)W^{u_1(k),\lambda_{Ev_1}}_{\varrho, \nu_1}(m) \overline{\chi_{w\nu}(m)W^{Au(k),\lambda_{EAv}}_{w\sigma, w\nu}(m)}\,dm\,dk\, d\nu\, d\mu_{\mathcal{P}}(\nu_1).\notag
\end{align}

As $\varrho$ is not equivalent to $w\sigma$ for any $w \in W(G|M_{\theta})$, the integral over $N_* \setminus M_{\theta}$ vanishes. 

Now we come to the final assertion of the theorem. In this we will assume that $P_{\theta}=P_{\vartheta}$ and because of Lemma \ref{lemma621},  $\mathcal{P}=\mathcal{O}$. However we will consider the function $W_{[\beta',\mathcal{O}, \bar{P}_{\theta}];v}$ where $\beta'$ is replaced by $$\beta' = (\nu \mapsto \beta(\nu) \lambda)$$  where $\beta \in C^{\infty}(\mathcal{O})$ and $\lambda \in Wh_{\psi^{\theta}}(H_{\sigma})$.

Now from \eqref{reff5} and \eqref{reff14} \begin{align*}\la W_{[\beta',\mathcal{O}, \bar{P}_{\theta}];u_1} ,& W_{[\Lambda(\alpha \otimes v), \mathcal{O}, \bar{P}_{\theta}];u} \ra = \\
&\int_{N_0 \setminus G}\int_{\mathcal{O}} \overline{J_{\sigma, \nu}(\beta(\nu)\lambda)(\pi_{\bar{P}_{\theta},\sigma, \nu}(g)u_1)}(f_{[\alpha, \mathcal{O}, \bar{P}_{\theta}];u,v})_{\psi}(g)\,\dsig\,dg\\
&=\int_{N_2 \setminus M_{\theta}\times K} \int_{\mathcal{O}}  \overline{\beta(\nu)\lambda(\chi_{\nu}\boxtimes\sigma(m)u_1(k))} ((R(k)(f_{[\alpha, \mathcal{O}, \bar{P}_{\theta}];u,v})_{\psi^{\theta}})^{P_{\theta}}(m)\,\dsig\, dm\,dk \\
&=\gamma(G|M_{\theta})c(G|M_{\theta})^2\int_{\mathcal{O}} \overline{\beta(\nu)}\int_{\mathcal{O}}\alpha(\nu_1) \\
&\qquad\int_{N_2 \setminus M_{\theta} \times K} \sum_{w \in W(G|M_{\theta})} \overline{\chi^{}_{w\nu}(m)W^{Au_1(k), M^{-1}\lambda}_{w\sigma, w\nu}(m)}\chi_{w\nu_1}(m)W^{Au(k), \lambda_{EAv}}_{w\sigma, w\nu_1}(m)\,dm\,dk\, d\nu_1\dsig
\intertext{by \eqref{reff12} and \eqref{reff14}} 
\end{align*}

The integral over $N_2 \setminus M_{\theta}$ vanishes unless $\nu=\nu_1$ in which case we are simply integrating on the image of the diagonal embedding $$ \mathcal{O} \hookrightarrow \mathcal{O} \times \mathcal{O}.$$ It is clear that the measure on the image is   $\dsig$. Then the expression above reduces to 
\begin{align}
 \gamma(G|M_{\theta})c(G|M_{\theta})^2\sum_{w\in W(G|M_{\theta})}\int_{\mathcal{O}}\alpha(\nu)\overline{\beta(\nu)} \int_{N_2 \setminus M_{\theta} \times K} \overline{W^{Au_1(k), M^{-1}\lambda}_{w\sigma, w\nu}(m)}W^{Au(k), \lambda_{EAv}}_{w\sigma, w\nu}(m)\, dm\,dk\, \dsig. \label{reff15}
\end{align}
Let $v' \in H_{w\sigma}$ be chosen to satisfy $ M_w(\nu)^{-1}\lambda(v')=1$ and set $\eta = \lambda_{v'}$. Then we have that $$W^{Au(k), \lambda_{EAv}}_{w\sigma,w\nu}(m)=\overline{\lambda_{w\sigma(m)Au(k)}(EAv) }=\overline{J_{w\sigma,w\nu}(\lambda_{w\sigma(m)Au(k)})(Av)}$$
and
$$\lambda_{w\sigma(m)Au(k)}=\overline{\eta(w\sigma(m)Au(k))}M_w(\nu)^{-1}\lambda$$
so that $$W^{Au(k), \lambda_{EAv}}_{w\sigma,w\nu}(m)=\eta(w\sigma(m)Au(k))\overline{J_{w\sigma,w\nu}(M_w(\nu)^{-1}\lambda)(Av)}.$$ Now,
\begin{align*}
\int_{N_2 \setminus M_{\theta} \times K} &\overline{W^{Au_1(k), M^{-1}\lambda}_{w\sigma, w\nu}(m)}W^{Au(k), \lambda_{EAv}}_{w\sigma, w\nu_1}(m)\,dm\,dk \\ &=\overline{J_{w\sigma,w\nu}(M_w(\nu)^{-1}\lambda)(A_w(\nu)v)}(\eta,M_w(\nu)^{-1}\lambda)_{w\sigma}\int_{ K}\la A_w(\nu)u(k),A_w(\nu)u_1(k)\ra_{w\sigma} \\
&=d(w\sigma)^{-1}\overline{J_{w\sigma,w\nu}(M_w(\nu)^{-1}\lambda)(A_w(\nu)v)}\la u,u_1\ra_{I(\sigma)}\\ \intertext{by Lemma \ref{bilinear identity} and since $\,^{\circ}c_{\bar{P}|\bar{P}}(w,\sigma)$ is unitary. Continuing from \eqref{reff15}, this expression simplifies to } 
\gamma(G|M_{\theta})&c(G|M_{\theta})^2\la u,u_1\ra_{I(\sigma)}\sum_{w\in W(G|M_{\theta})}d(w\sigma)^{-1}\int_{\mathcal{O}}\alpha(\nu)\overline{\beta(\nu)J_{w\sigma,w\nu}(M_w(\nu)^{-1}\lambda)(A_w(\nu)v)} \,\dsig\\
\intertext{and by \eqref{reff14} and the well known fact that $d(w\sigma)=d(\sigma)$, this equals}
&\gamma(G|M_{\theta})c(G|M_{\theta})^2\la u,u_1\ra_{I(\sigma)} |W(G|M_{\theta})|\int_{\mathcal{O}} d(\sigma)^{-1}\alpha(\nu)\overline{ \beta(\nu)J_{\sigma,\nu}(\lambda)(v)} \,\dsig\\
&=\gamma(G|M_{\theta})c(G|M_{\theta})^2\la u,u_1\ra_{I(\sigma)} |W(G|M_{\theta})|\int_{\mathcal{O}}(\alpha'(\nu),\beta'(\nu))_{\sigma} \,\dsig
\end{align*}
where $\alpha'(\nu) = \Lambda(\alpha \otimes v)$ and $\beta'(\nu)= \beta(\nu)\lambda$ for $\alpha, \beta \in C^{\infty}(\mathcal{O})$. Lemma \ref{lastindentity} justifies the last equality. 

The proof is now complete.
\end{proof}


\subsection{}
Fix an open compact subgroup $H_m$, $f \in \mathcal{C}^*(N_0 \setminus G; \psi)^{H_m}$ and $\nu \in \mathrm{Im}X_{ur}(M_{\theta})$.  Observe that the integral in Corollary \ref{pairingconv} defines a conjugate linear functional on $$(I^G_{\bar{P}_{\theta}}\sigma \boxtimes \chi_{\nu})^{H_m}\otimes Wh_{\psi^{\theta}}(H_{\sigma})=:\mathscr{H}^{H_m}_{\sigma,\nu}$$ considered as a Hilbert space.  By the Reisz Representation theorem there exists a unique element in $I(\sigma) \otimes Wh_{\psi^{\theta}}(H_{\sigma})$ denoted $W_{\bar{P}_{\theta}, \sigma}(f)(\nu)$ such that 
$$ \la W_{\bar{P}_{\theta}, \sigma}(f)(\nu), w\otimes \lambda \ra_{\mathscr{H}^{H_m}_{\sigma,\nu}}=\int_{N_0 \setminus G} f(g) \overline{J_{\sigma}(\lambda)(\pi_{\bar{P}_{\theta}, \sigma, \nu}(g)w)}\, dg.$$

As we vary $\sigma\boxtimes \chi_{\nu}$ over $\mathcal{O}$, we may consider $W_{\bar{P}_{\theta}, \sigma}(f)$ as an element of $I(\sigma) \otimes C^{\infty}(\mathcal{O}; Wh_{\psi^{\theta}}(H_{\sigma}))$.  

Define the following linear operator on $I(\sigma) \otimes C^{\infty}(\mathcal{O}; Wh_{\psi^{\theta}}(H_{\sigma}))$ to $\mathcal{C}^*(N_0 \setminus G;\psi)$, $$T_{\bar{P}_{\theta}, \sigma}(v \otimes \alpha):= |W(G|M_{\theta})|^{-1} \gamma(G|M_{\theta})^{-1}c(G|M_{\theta})^{-2}W_{[\alpha, \mathcal{O}, \bar{P}_{\theta}]; v}.$$
Let $\mathscr{E}_{\psi^{\theta}}^2(M_{\theta})$ denote the set of isomorphism classes of $\psi^{\theta}$ generic discrete series representations on $M_{\theta}$. By $\mathscr{E}_{\psi^{\theta}}^2(M_{\theta})/W(G|M_{\theta})$ we mean the quotient space obtained by identifying the discrete series representations which are conjugate under the natural action of the Weyl group $W(G|M_{\theta})$. 

\begin{thm}[Plancherel Formula for $L^2(N_0 \setminus G;\psi)$]
Consider $f \in \mathcal{C}^*(N_0 \setminus G; \psi)$, then\begin{align} f = \sum_{(P_{\theta},A_{\theta}) \succ (P_0 , A_0)} \sum_{\sigma \in \mathscr{E}_{\psi^{\theta}}^2(M_{\theta})/W(G|M_{\theta})} T_{\bar{P}_{\theta}, \sigma}(W_{\bar{P}_{\theta}, \sigma}(f)).\label{plan}\end{align}
\end{thm}
\begin{proof}Set $\mathscr{H}_{\sigma,\nu} = \bigcup_{K} \mathscr{H}^K_{\sigma, \nu}$ and consider the Hilbert space defined by $$\mathcal{I}_{\sigma,M_{\theta}}=\int^{\oplus}_{\mathcal{O}} \mathscr{H}_{\sigma, \nu}\, d\tilde{\mu}_{\mathcal{O}}(\nu)$$ where $$d\tilde{\mu}_{\mathcal{O}}= |W(G|M_{\theta})|^{-1} \gamma(G|M_{\theta})^{-1}c(G|M_{\theta})^{-2} \dsig.$$ 

Then we consider $I(\sigma)\otimes C^{\infty}(\mathcal{O};Wh_{\psi^{\theta}}(H_{\sigma}))$ as a dense subspace of $\mathcal{I}_{\sigma}$. The linear operator $T_{\bar{P}_{\theta}, \sigma}$  extends to $\mathcal{I}_{\sigma}$ by continuity. 

We check that for any $v_1 \otimes \alpha_1, v_2\otimes \alpha_2 \in I(\sigma)\otimes C^{\infty}(\mathcal{O};Wh_{\psi^{\theta}}(H_{\sigma}))$, 
\begin{align*}
\la T_{\bar{P}_{\theta}, \sigma}(v_1\otimes \alpha_1), T_{\bar{P}_{\theta}, \sigma}&(v_2 \otimes \alpha_2) \ra_{L^2} \\
&= |W(G|M_{\theta})|^{-2} \gamma(G|M_{\theta})^{-2}c(G|M_{\theta})^{-4}\la W_{[\alpha, \mathcal{O}, \bar{P}_{\theta}]; v_1}, W_{[\alpha, \mathcal{O}, \bar{P}_{\theta}];v_2}\ra \\ 
\intertext{by definition of $T_{\bar{P}_{\theta}, \sigma}$} 
&=|W(G|M_{\theta})|^{-1}\gamma(G|M_{\theta})^{-1} c(G|M_{\theta})^{-2} \left(\int_{\mathcal{O}} (\alpha_1(\nu), \alpha_2(\nu))_{\sigma}\, \dsig\right)\la v_1, v_2 \ra_{I(\sigma)} \\ 
\intertext{by Theorem \ref{main4}(3)}
&=\left(\int_{\mathcal{O}} (\alpha_1(\nu),\alpha_2(\nu))_{\sigma}\, d\tilde{\mu}_{\mathcal{O}}(\nu)\right) \la v_1,v_2 \ra_{I(\sigma)} \\
&= \la v_1 \otimes \alpha_1 , v_2 \otimes \alpha_2 \ra_{\mathcal{I}_{\sigma,M_{\theta}}}.\end{align*}

Thus we see that $T_{\bar{P}_{\theta}, \sigma}$ extends to a unitary operator from $\mathcal{I}_{\sigma}$ into $L^2(N_0\setminus G; \psi)$. 

By Theorem \ref{main4} (2) and Lemma \ref{lemma621}, one sees that $$T:=\sum_{(P_{\theta},A_{\theta})\succ (P_0, A_0)} \sum_{\sigma \in \mathscr{E}_{\psi^{\theta}}^2(M_{\theta})/W(G|M_{\theta})} T_{\bar{P}_{\theta}, \sigma}$$ is a linear map from $$\sum_{(P_{\theta},A_{\theta})\succ (P_0, A_0)} \sum_{\sigma \in \mathscr{E}_{\psi^{\theta}}^2(M_{\theta})/W(G|M_{\theta})} I(\sigma) \otimes C^{\infty}(\mathcal{O}; Wh_{\psi^{\theta}}(H_{\sigma}))$$ \emph{onto} $\mathcal{C}^*(N_0 \setminus G;\psi)$. By the previous paragraph, $T$ is extends to \emph{unitary} linear operator on  $$\sum_{(P_{\theta},A_{\theta})\succ (P_0, A_0)} \sum_{\sigma \in \mathscr{E}_{\psi^{\theta}}^2(M_{\theta})/W(G|M_{\theta})} \mathcal{I}_{\sigma,M_{\theta}}.$$

As before, let $\alpha(\nu)=\alpha'(\nu)\lambda$ for any $\alpha \in C^{\infty}(\mathcal{O})$ and a fixed $\lambda \in Wh_{\psi^{\theta}}(H_{\sigma})$.  Then, we see that \begin{align*} 
\la T_{\bar{P}, \sigma}(W_{\bar{P}, \sigma}(f)), T_{\bar{P}, \sigma}( w \otimes \alpha) \ra_{L^2} &= \int_{\mathcal{O}}\overline{ \alpha'(\nu)} \la W_{\bar{P},\sigma}(f), w \otimes \lambda \ra_{\mathscr{H}_{\sigma, \nu}} d\tilde{\mu}_{\mathcal{O}}(\nu) \\
&= \int_{\mathcal{O}} \overline{\alpha'(\nu)} \int_{N_0\setminus G} f(g) \overline{J_{\sigma, \nu}(\lambda)(\pi_{\bar{P}_{\theta}, \sigma, \nu}(g)w)}\, dg \, d\tilde{\mu}_{\mathcal{O}}(\nu) \\
&=\la f, T_{\bar{P}, \sigma}(w\otimes \alpha)\ra_{L^2}.
\end{align*} 

It is clear that if $P \neq P'$ or if $\sigma$ and $\sigma'$ are not Weyl group conjugates of each other, then $\la T_{\bar{P}, \sigma}(W_{\bar{P}, \sigma}(f)), T_{\bar{P}', \sigma'}( w \otimes \alpha) \ra_{L^2}=0.$

Let $R(f)$ denote the right hand side of \eqref{plan}. Then the considerations of the previous two paragraphs imply that $\la f, h \ra = \la R(f), h\ra$ for all $h \in \mathcal{C}^*(N_0 \setminus G; \psi)$. As $\mathcal{C}^*(N_0 \setminus G; \psi)$ is dense in $L^2(N_0 \setminus G; \psi)$, $f= R(f)$ and the theorem follows. 
\end{proof}
We state a corollary of the proof. 
\begin{tcor}\label{lastcor}
The unitary linear map $T$ is a surjection from $$\sum_{(P_{\theta},A_{\theta})\succ (P_0, A_0)} \sum_{\sigma \in \mathscr{E}_{\psi^{\theta}}^2(M_{\theta})/W(G|M_{\theta})} \mathcal{I}_{\sigma,M_{\theta}}$$ onto $L^2(N_0 \setminus G; \psi)$. 
\begin{flushright}
$\qed$ 
\end{flushright}

\end{tcor}

\end{raggedright}


\end{document}